\documentclass[11pt]{article}
\usepackage{cleveref}    
\usepackage{amsmath}
\usepackage{amsfonts}
\usepackage{amsthm}
\usepackage{centernot}
\usepackage{eurosym}
\usepackage{graphicx}
\usepackage{amssymb}
\usepackage{epstopdf}
\usepackage{topcapt}
\usepackage{latexsym}
\usepackage{mathtools}
\usepackage{graphics}
\usepackage{subfigure}
\usepackage{dsfont}
\usepackage[english]{babel}
\long\def\ABSTRACT#1 {\parindent=0pt {\rm #1 }\parindent=1.5cm}

\newtheorem{theorem}{Theorem}[section]

\newtheorem{corol}[theorem]{Corollary}
\newtheorem{prop}[theorem]{Proposition}

\newtheorem{remark}[theorem]{Remark}

\newtheorem{hypo}[theorem]{Hypothesis}

\newtheorem*{details*}{Details}

\newtheorem{lemma}[theorem]{Lemma}
\newtheorem*{lemma*}{Lemma}
\newtheorem{definition}[theorem]{Definition}

\newtheorem*{reminder*}{Reminder}

\def\sqr#1#2{{\vcenter{\vbox{\hrule height.#2pt
\hbox{\vrule width.#2pt height#1pt \kern#1pt
\vrule width.#2pt}
\hrule height.#2pt}}}}

\def\a{\alpha}

\def\de{\delta}
\def\e{\varepsilon}

\def\g{\gamma}
\def\8{\infty}
\def\k{\kappa}

\def\nn{\nonumber}

\def\ome{\omega}
\def\Ome{\Omega}

\def\ot{{\overline \tau}}

\def\obeta_2{{\overline \beta_2}}

\def\ov{{\overline v}}

\def\ox{{\overline x}}

\def\0*{^{\odot *}}
\def\qed{\hfill$\sqr45$\bigskip}

\def\N{{\mathds N}}
\def\R{{\mathds R}}

\def\r{\rho}
\def\s{\sigma}

\def\t{\tau}
\def\th{\theta}
\def\vi{\varphi}

\def\Proof{\noindent{\em Proof. }}
\def\<{\langle}
\def\>{\rangle}
\def\->{\longrightarrow}
\def\_>{\rightarrow}

\title{Global extinction, dissipativity and persistence for
a certain class of differential equations with state-dependent delay}
\author{Philipp Getto \and Gergely R\"ost}
\date{}
\begin{document}
\maketitle


\tableofcontents

%
\section{Introduction}
%
%
%
%
\subsection{Overall summary} 
In this paper we study, at different levels of generality, certain systems of delay differential equations (DDE). One focus and motivation is a system with state-dependent delay (SD-DDE) that has been formulated to describe the maturation of stem cells. We refer to this system as the cell SD-DDE. In the cell SD-DDE, the delay is implicitly defined by a threshold condition. The latter is specified by the time at which the (also implicitly defined) solution of an external nonlinear ordinary differential equation (ODE), which is parametrised by a component of the SD-DDE, meets a given threshold value. We investigate the dynamical properties global asymptotic stability (GAS) of the zero equilibrium, persistence and dissipativity/ultimate boundedness. 
\subsection{The equations under study}
The cell SD-DDE takes the form
%
\begin{eqnarray}
w'(t)&=&q((v(t))w(t), 
\label{eq68}\\
v'(t)&=&
\frac{\g(v(t-\t(v_t)))g(x_2,v(t))}{g(x_1,v(t-\t(v_t)))}e^{\int_0^{\t(v_t)}(d-D_1g)(y(s,v_t),v(t-s))ds}w(t-\t(v_t))
\nonumber\\
&&-\mu v(t),
\label{eq11}
\end{eqnarray}
for $t\ge 0$ and 
\begin{eqnarray}
(w,v)_0&=&(\vi,\psi).
\label{eq70}
\end{eqnarray}
\noindent
We use the standard notation
\begin{eqnarray}
x_t(s):=x(t+s),\;\;s<0,
\label{eq73}
\end{eqnarray}
if a function $x$ is defined in $t+s\in{\bf R}$. If $t$ is fixed, then $x_t$ is a function describing the history 
of $x$ at time $t$. We use subindices also in other contexts than histories and hope that our verbal
explanations will avoid confusion. 
Both, (\ref{eq68}) and (\ref{eq11}), are equations in 
${\bf R}$ and all functions are real-valued. The functions $q$, $\g$, $g$ and $d$ have real arguments, $
\mu$ is a parameter and $\g$, $g$, $d$ (the second $d$ in the exponent in (\ref{eq73}) refers of course to integration), $\t$ and $\mu$ take nonnegative values. We use $D_i$ to 
denote the derivative of a function with respect to the $i^{th}$ argument. The delay $\t$ depends on the second component $v_t$ of the state (only). It is implicitly
given by a nonlinear functional: For a function $\psi:[-h,0]\->\R$, $h>0$, we define $\t=\t(\psi)\in[0,h]$ as the solution of the equation 
\begin{equation}
y(\t,\psi)=x_1,
\label{eq7}
\end{equation}
where $y(\cdot,\psi)$ is defined via the ODE
\begin{eqnarray}
\begin{array}{ccc}
y'(s)&=&-g(y(s),\psi(-s)),\;\;s>0,
\\
y(0)&=&x_2,
\end{array}
\label{eq5}
\end{eqnarray}
and $x_1, x_2\in{\bf R}$, $x_1<x_2$ are given model parameters (and hence subindices are used
in another context than in (\ref{eq73})). 
We define $y$ going backward in time to facilitate denoting time dependence in the second argument of 
$g$, given that $\psi$ is defined on 
$[-h,0]$.

%
%
%
The system describes the dynamics of a stem cell population of size $w(t)$, regulated by the mature cell 
population of size $v(t)$. The equations have been 
deduced  via integration along the characteristics from a partial differential equation describing the 
``transport" of a density $n(t,\xi)$, where $\xi\in[x_1,x_2]$ is the maturity of committed cells (former stem cells
that have started the maturation process, but are not fully mature). See \cite{Getto}
and references therein for the latter and \cite{Getto1, Vivanco, Scarabel} for modelling and biological 
background. 
\subsection{Own previous results  cell SD-DDE} 
%
%
\subsubsection{Terminology}
For precise definition of some terminology used in the remainder of the introduction we refer to Section \ref{ss3}. An equilibrium is called 
{\it globally asymptotically stable}, if it is globally attractive and stable. We write GAS for abbreviation, also 
for related word combinations like {\it global asymptotic stability}. 
\subsubsection{ODE model in MB} 
In the cell SD-DDE the limit case $x_1=x_2$ corresponds to $\t\equiv 0$, hence the term in the exponent vanishes too and the DDE
becomes a two-dimensional ODE. This ODE, which can be interpreted as a simplistic model on its own, 
has been studied in \cite{Nakata}. There is always a zero equilibrium. In a region of the parameter space, there also exists a unique positive equilibrium, in a disjoint region there is no positive equilibrium.
The region in which the positive equilibrium exists is often called {\it existence  region}.
It has been shown that zero is GAS in 
absence of the positive equilibrium and unstable in the existence region. Moreover, a Lyapunov function 
has been constructed to prove that the positive equilibrium is GAS throughout the existence region. 
%
\subsubsection{WP on manifold for cell SD-DDE} 
The cell SD-DDE is of the form
\begin{eqnarray}
x'(t)=F(x_t),\;t>0, 
\label{eq55}
\end{eqnarray}
where $F$ has two-dimensional range space. If $F$ has $N$-dimensional range space, we will write
about $N$-dimensional DDE.  For the case, where $F$ fulfils certain smoothness conditions on an open
subset of $C$ (see the opening paragraph of Section \ref{ss4} for notation and related definitions), the classical
theory of functional differential equations (FDE) of \cite{HaleLunel} applies. For SD-DDE, it is known that the
hinted conditions are too strong, see e.g. \cite{Walther1}. On the other hand, if $F$ satisfies 
certain smoothness conditions on an open subset of $C^1$, well-posedness for $N$-dimensional DDE on the solution manifold, a
continuously differentiable sub-manifold of $C^1$, along with the fact that solutions define a maximal 
differentiable semiflow, is established in \cite{Walther1}. The approach is tailored to 
achieve differentiability for semiflows induced by SD-DDE. If we here write about general DDE, we mean 
this setting. 

In \cite{Getto}, the smoothness conditions are guaranteed for the functional inducing 
the cell SD-DDE and, by application of the general result, well-posedness and the existence of a maximal 
differentiable semiflow has been established for the cell SD-DDE. In the same paper a criterion for global 
existence for general DDE is established and applied to the cell SD-DDE, such that for the latter also 
associability of a global differentiable semiflow is proven. We will reformulate the result below. 
%
%
\subsubsection{Equilibria and stability for cell SD-DDE} 
Equilibria of the cell SD-DDE and their local stability has been studied analytically and numerically in 
\cite{Scarabel}. Similar as in the two-dimensional ODE system, there always is a trivial 
equilibrium and for the unique positive equilibrium there is an existence region of the parameter space 
and a region of non-existence. The trivial equilibrium is locally stable in the region in which the positive 
does not exist and destabilises  in a transcritical bifurcation upon entering the existence region. Close to 
the bifurcation point, the positive equilibrium is locally stable. However, away from the bifurcation point, other than in the ODE system, it may destabilise (evidentially in a Hopf bifurcation). The analytical stability 
results are based on linearised stability theorems for the cell SD-DDE. These have been established by 
combining the differentiability of the functional proven in \cite{Getto} with the linearised stability 
theorems for general DDE in \cite{Hartung} (stability) and \cite{Stumpf} (instability).  
%
%
\subsection{Motivation for this research: GAS of zero in absence of positive, persistence in presence} 
The aim of this research is to show at different levels of generality that, in the existence region of the positive equilibrium, there is persistence and in the region of non-existence the stability of the zero steady state is global. 
%
%
\subsubsection{Existing persistence methods} 
For persistence theory we will refer to the monograph \cite{Smith}. In the preface,  the notion of persistence is introduced to describe the survival 
of a species over a long term. In the interpretation of our equations this suggests to replace the notion
of species by the cell-types that $w$ and $v$ describe. One result of the present paper, uniform persistence with respect to the persistence function 
\begin{eqnarray}
\r_m(\vi,\psi):=\min\{\vi(0),\psi(0)\},
\label{eq10}
\end{eqnarray}
can be formulated as the existence of some $\e>0$, such that 
\[
\liminf_{t\->\infty}\min\{w(t),v(t)\}>\e.
\]
If the $\liminf$ is replaced by a $
\limsup$, the resulting property corresponds to uniform weak persistence of our system. Since the 
persistence theory in \cite{Smith} is for semiflows on metric spaces and the solution manifold for SD-DDE 
induces a metric space in $C^1$, we should apply the theory to the $C^1$-topology. An important 
theoretical result in \cite{Smith} is that uniform persistence can be concluded essentially 
from uniform weak persistence and a compactness property. The latter could be described as dissipativity/ultimate 
boundedness, where the bounded set additionally is compact, or as the existence of a compact pointwise attractor 
of the state space. A further result in \cite{Smith} allows to conclude weak persistence under the assumption of a 
very related compactness property. It turned out that the structure of the cell SD-DDE brings some challenge to the 
establishment of these compactness properties. 
%
%
\subsubsection{Existing dissipativity methods and why they don't work here} 
Recall, that in the limit case of a two dimensional ODE it has been established that 
there is always a GAS equilibrium, such that ultimate boundedness follows trivially.  On the 
other hand, we discussed that for the (general) cell SD-DDE the positive  equilibrium may be unstable. 
This hints at problems with a generalisation of the construction of a Lyapunov function. 

A class of systems, for which methods to show dissipativity are well-known are certain systems describing 
consumer-resource or chemostat like behaviour. There is one species that in absence of others has a self-regulatory mechanism, whereas the presence of others has only negative influence. The other species
depend on the former, e.g. via consumption 
The idea is  then to first show that the state-component associated with the self-regulatory species is ultimately 
bounded and then deduce ultimate boundedness of the  components corresponding to the other species. 
The underlying theoretical result roughly says that dominance of a (component of a) solution can be 
concluded from a corresponding dominance of right hand sides. For ODE such a result can be found in 
\cite{Hale}. 
For DDE such results are known if additional monotonicity properties can be referred to \cite[Theorem 4.1]{Hirsch}.

Another method that can sometimes be applied 
is related to conservation of mass 
considerations. It is first shown that the sum of the components is ultimately bounded 
and then conclusions for the components are drawn.

In \cite[4.1.3]{Hale} dissipativity is shown for a  one-component DDE that was already studied in \cite{Wright}. In \cite[4.1.3]{Hale}, it is distinguished between the case where the set of time points at which equilibrium is reached is bounded and the case, where it is unbounded. 
When trying to adapt the proof to some special cases of the cell SD-DDE (fixed delay), the authors of the present paper so far did not manage to resolve the case where a similar set is unbounded. 

In summary, for the authors, even for the case of a fixed delay, none of the above methods worked for the cell 
SD-DDE. A hint at the problems is the combination of the existence of the delay with the specific coupling of the components, which features one component, i.e. $v$,  that typically inhibits both components and another component, i.e., $w$, that stimulates both.  
%
%
\subsection{Main results and methodology} 
%
%
\subsubsection{Introduction, motivation and specification of object 2, the cell FDE} 
We have found, that the methods that we elaborate here essentially apply to a certain class of 
two-dimensional DDE, that is more general than the cell SD-DDE and that can be formulated as (\ref{eq68}) coupled to
\begin{eqnarray}
v'(t)=j(w_t,v_t)-\mu v(t),  
\label{eq58}
\end{eqnarray}
for a general functional $j$. We call (\ref{eq68}, \ref{eq58}) the {\it cell FDE}. We find it instructive to present the methods for both, 
the cell FDE and the cell SD-DDE, in a top down approach. In the methods presented for the cell FDE we often 
use the hierarchical structure given by the two components and the linearity of (\ref{eq68}), 
such that at the corresponding places, we did not aspire more generality at this point. Having specified
two of the mathematical objects under study, in the following we outline our results and methods. 
%
%
\subsubsection{Main result 1 and methodology GAS of zero} 
Our first main result is the elaboration of general conditions for GAS of zero in the $C$-norm and in the 
$C^1$-norm for a certain region of the parameter space. 

To achieve this, we first give some conditions for general DDE under which notions of stability on a given
subset of the state space in the $C$-topology imply the corresponding notions in the $C^1$-topology.

We then 
focus on the cell FDE and the case where $j(0)=0$, i.e., where the cell FDE has a a zero solution. This case applies to the cell SD-DDE. Then, we show  
that a 
continuity property of $j$ -- in the $C$-topology -- in zero  can be used to transfer notions of stability and 
attractivity of the zero equilibrium from the $w$-component to the system. Since the stability properties for the $w$-component are assumed, for the proofs only the structure of the second equation of the cell FDE is relevant. 
The last property, i.e., transferability of  stability notions from $w$-component to cell FDE, will be used to establish both, GAS of zero and persistence results. 

Finally,  
 we focus on the case where $q$ satisfies a certain negativity property. The property
implies that $q$ is negative for positive arguments and nonpositive in zero. 
Then, linearity of the $w$-equation allows to establish stability,  and if $q(0)<0$ also attractivity, of zero for 
the $w$-component rather quickly. If $q(0)=0$, we will show that attractivity is not always given. 
Under a general growth condition for $j$ 
however, we can guarantee attractivity for the $w$-component, which completes the GAS of zero analysis. 
%
%
\subsubsection{Main result2 and methodology compact attractor} 
Our second main result is the elaboration of a method to establish above discussed compactness properties. The method is tailored to the cell FDE. We first establish some auxiliary results for general DDE: A Lipschitz 
property of the functional can be used to conclude a compactness property from dissipativity in the $C^1$-norm. 
Moreover, a boundedness property of the functional can be 
used to transfer dissipativity from the $C$-norm to the $C^1$-norm. 
Then, for the cell FDE, it turned out convenient to establish compactness properties by showing dissipativity in 
both norms directly, i.e., without the boundedness property for the functional. For the cell 
SD-DDE, we can establish dissipativity in both norms directly or combine the boundedness property of the 
functional with dissipativity in the $C$-norm. 

Then we focus on the cell FDE. First, we assume that $q$ satisfies a property of ultimate negativity
(which relaxes the assumption for $q$ made in the GAS  of zero analysis). Then we suppose that 
$j$ can become arbitrarily large, independent of the $v$-argument, if the $w$-argument is sufficiently large. As a result, very roughly speaking, if $w$ grows, then $v$ grows, and $w$ goes back to decreasing. A precise elaboration is somewhat subtle. We outline it in the following. 

We first show that there is a constant such that, if the $w$-component remains long enough above it, it 
starts to decrease. Under the above assumption on $j$, constant and decrease duration are uniform with 
respect to initial conditions.  On the other hand, the $w$-component obviously cannot increase faster than 
exponentially. Hence, if another constant is chosen so large that with exponential increase it could be 
reached only in a time longer than the decrease duration, the second constant can never be reached 
and dissipativity of the $w$-component in the $C$-norm follows. 

To guarantee global existence of solutions,  we assume that $j$ maps, in the first argument  uniformly with 
respect to the second, bounded sets 
on bounded sets. Note that, since for the cell SD-DDE  $j$ 
is linear in $w$, boundedness of $j$ in general is unclear. With the same boundedness assumption for 
$j$, we can conclude dissipativity of the $v$-component, hence of the tuple, in the $C$-norm. 
Then the differential equation can be used to transfer dissipativity to the $C^1$-norm. Compactness properties
can be concluded. 
%
%
\subsubsection{Main Result 3  and methodology Persistence} 
Our third and final main result is the elaboration of conditions for persistence, 
essentially in the region of the parameter space complementary to the region in which zero is GAS. 
If discussed compactness is provided, uniform persistence can be guaranteed, if uniform weak 
persistence is established. We have mentioned that the case where $q(0)<0$ and $q$ satisfies negativity
properties can be associated with GAS of zero. In this sense, we cannot expect 
persistence for this case. 
Towards a persistence proof, we first assume that $q(0)>0$. 
On the set defined by $\vi(0)=0$, in the notation of (\ref{eq70}), zero is 
stable in the $C$-norm. As discussed, this implies stability of zero in $C^1$. If $\vi(0)>0$, we can show 
that $q(0)>0$ always leads to the impossibility of convergence to zero. In \cite{Smith} weak repellence 
and uniform weak persistence are defined in terms of persistence functions and repellence of zero is prescribed as a precondition for uniform weak persistence. This motivates finding persistence functions $\r$ such that $\vi(0)>0$ can be related to weak $\r$-repellence of zero and uniform 
weak $\r$-persistence. 
We show weak persistence for the persistence function 
\begin{eqnarray}
\r_1(\vi,\psi)=\vi(0)
\label{eq16}
\end{eqnarray}
via an elementary contradiction proof. To show weak persistence for the persistence function $\r_m$ 
we use the theory of \cite{Smith}, which requires first establishing compactness properties.

In the analysis of the three properties, GAS of zero, dissipativity and persistence, we typically 
encountered cases where the property holds for both components. On the other hand, again for the three, 
the methods turned out to first focus on the $w$-component and then establish inheritance for the $v$-component. 

%
%
\subsection{Short sectionwise summary} 
For better readability, we tried to give preference to a top-down structure in the following way,
even though this implied starting with a subsection, which for the majority consists of 
established results: In Subsection \ref{ss3} we summarise established notions and theory of persistence, 
dissipativity and stability in a setting of given semiflows on metric spaces. In Subsections \ref{ss4} and 
\ref{ss2} we elaborate some preliminary results for general DDE and for the cell FDE, respectively. In 
particular, we guarantee that the systems define semiflows.  We remain with the cell FDE in 
Sections \ref{s2} (GAS of zero) and \ref{s4} (persistence and dissipativity). Then, in Section \ref{ss1}
we focus on the cell SD-DDE and Section \ref{s6} closes with examples and discussion. 
%
\subsection{General remarks}
Unless stated otherwise, a ``Hypothesis" or an assumption in the running text is assumed to hold from the point where it is formulated to the end of the subsection or, if there is no subsection structure, to the end of the section and will be tacitly dismissed
beyond. To enable quicker reference, in ``Theorems" all assumed ``Hypotheses" will be referred to explicitly. 

%
%
\section{From semiflows on metric spaces to the cell FDE}
%
%
\subsection{Semiflows on metric spaces - definitions and established results}\label{ss3}
For the following notions, definitions and related information, see \cite[Definitions 2.1, 2.4, 2.25 and 8.15] {Smith}. Let $Y$ be a metric space with metric $d$,
$\R_+:=[0,\infty)$ and $\Phi:\R_+\times Y\->Y$ a continuous semiflow. We define the distance between a point
$x$ and a set $B$ and the distance between two sets $A$ and $B$, respectively, as
\begin{eqnarray}
d(x,B):=\inf\{d(x,y):\;y\in B\},\;d(A,B):=\sup\{d(x,B):\;x\in A\}. 
\nn
\end{eqnarray}

By a {\it persistence function} we mean an arbitrary function $\rho:Y\->\R_+$ that is continuous and non-identically zero. We will use $\r$ and $\tilde\r$ to denote persistence functions.
\begin{definition}\label{d1}\rm
\begin{itemize}
\item[(i)] Let $A\subset Y$ and let $\{A_t\}_{t\in\R_+}$ be a family of subsets of $Y$. We say that $A_t$ {\it converges} to $A$, symbolically $A_t\->A$, if for any open set $U\subset Y$, $A \subset U$, there exists some $r=r(U)\ge 0$, such that $A_t\subset U$
for all $t\in[r,\infty)$. 
\item[(ii)] A set $K\subset Y$ is said to {\it attract} a set $M\subset Y$, if $K\neq\emptyset$ and 
\[
d(\Phi(t,M),K)\->0,\;\;{\rm as}\;t\_>\infty.
\]

\item[(iii)] $\Phi$ is called {\it point dissipative} or {\it ultimately bounded} if there exists a bounded subset
$B$ of $Y$, which attracts all points of $Y$. 

\item[(iv)] A set $M\subset Y$ is called
{\it weakly $\rho$-repelling}, if there is no $x\in Y$, such that $\rho(x)>0$ and $\Phi(t,x)\->M$
as $t\_>\infty$.
\end{itemize}
\end{definition}
If $A_t, A\neq\emptyset$ for all $t\ge 0$, then $A_t\->A$ implies $d(A_t,A)\->0$, 
whereas the converse is only true, if $A$ is compact \cite[Lemma 2.3]{Smith}. 
In this paper, we typically consider families of one-pointed sets of the form $\Phi(t,x)$. 
In order to apply a criterion for persistence, we define
\[
X_0^\r:=\{\phi\in Y:\;\rho(\Phi(t,\phi))=0\;\forall\;t\ge 0\},
\]
see \cite[Section 8.3]{Smith}. 
Note that, if $\r\ge\tilde\r$, then $X_0^\r\subset X_0^{\tilde\r}$.  We can formulate \cite[Definition 2.7]{Smith}
for a semiflow on $\R_+$ as
\begin{definition}\label{d2}\rm
The {\it $\ome$-limit set}
of a subset $M\subset Y$ is defined as
\[
\ome(M)=\bigcap_{t\ge 0}\overline{\Phi([t,\infty)\times M)}.
\]
\end{definition}
We adapt \cite[Definition 3.1]{Smith} to our setting:
\begin{definition}\rm  The semiflow $\Phi$ is called {\it uniformly weakly $\rho$-persistent} if there is some $\e>0$ such that 
\[
\limsup_{t\_>\infty}\rho(\Phi(t,x))> \e,\;{\rm whenever}\;x\in Y,\;\rho(x)>0.
\]
It is called {\it uniformly (strongly) $\r$-persistent} if there exists some $\e>0$, such that 
\[
\liminf_{t\_>\infty}\r(\Phi(t,y))>\e,\;{\rm whenever}\;y\in Y, \;\r(y)>0. 
\]
\end{definition}
The proof of the following result is straightforward and we omit it. 
\begin{lemma}\label{lem17}
Suppose that $\r(\phi)\ge \tilde\r(\phi)$ for all $\phi\in Y$, then $\Phi$ is uniformly weakly (strongly) $\r$-persistent, if it is  uniformly weakly (strongly) $\tilde\r$-persistent. If the ordering holds on a subset of $Y$, then the subset is weakly $\tilde\r$-repelling, if it is weakly $\r$-repelling. 
\end{lemma}
We will use \cite[Definitions 1.39 and 8.14]{Smith}
for semiflows on $\R_+$:
\begin{definition}\rm
(i) A map $T:\R\->Y$ is called a {\it total
trajectory} if $\Phi(t,T(s))=T(t+s)$ for all $t\ge 0$, $s\in\R$. 
\\
(ii) Let $C,B\subset X_0^\r$. $C$ is said to be {\it chained} to $B$ in $X_0^\r$, written $C\mapsto B$, if there is a total trajectory $T$ in $X_0^\r$ with $T(0)\notin C\cup B$, $T(-t)\->C$ and $T(t)\-> B$ as $t\_>\infty$. A finite collection $\{M_1,...,M_k\}$ of subsets of $X_0^\r$
is called {\it cyclic} if, after possible renumbering $M_1\mapsto M_1$, or $M_1\mapsto M_2\mapsto...\mapsto M_j\mapsto M_1$ in $X_0^\r$ for some $j\in\{2, ...., k\}$. Otherwise it is called {\it acyclic}. 
\end{definition}
We define
\[
\Ome:=\bigcup_{x\in X_0^\r}\ome(x),
\]
see \cite[(8.2)]{Smith}. 
We reformulate \cite[Theorem 4.5]{Smith}. 
\begin{theorem}\label{theo1}
Let $\s:=\r\circ\Phi:\R_+\times Y\->\R_+$. Assume that 
\\
(i) $\s$ is continuous. 
\\
Suppose that there exists a nonempty compact set $B\subset Y$ such that
\\
(ii) there are no $y\in B$, $s,t\in [0,\infty)$, such that $\r(y)>0$, $\s(s,y)=0$, 
$\s(s+t,y)>0$,
\\
(iii) for all $x\in Y$ one has $d(\Phi(t,x),B)\-> 0$ as $t\->\infty$. 
\\
Then the semiflow $\Phi$ is uniformly $\r$-persistent whenever it is uniformly weakly $\r$-persistent. 
\end{theorem}
The formulation of \cite[Theorem 4.5]{Smith}  is slightly stronger than in the previous theorem, in the way that in the former
(iii) is required to hold only for those $x$ with $\r(x)>0$. On the other hand, below we show the property for
all $x\in Y$. Moreover, if the result holds for all $x\in Y$, then, by our remark after Definition \ref{d1}, it trivially
implies (H) in \cite[Section 8.3]{Smith}. The latter is a precondition for \cite[Theorem 8.17]{Smith}, which we
will apply and reformulate as follows. 
\begin{theorem}\label{theo3}Suppose that Theorem \ref{theo1} (iii) holds for some nonempty compact $B\subset Y$, $\Ome\subset\bigcup_{i=1}^k M_i$, where 
each $M_i\subset X_0^\r$ is isolated (in $Y$), compact, invariant and weakly $\r$-repelling, $M_i\cap M_j=\emptyset$, if $i\neq j$. If $\{M_1,...,M_k\}$ is acyclic, then $\Phi$ is uniformly weakly 
$\r$-persistent. 
\end{theorem}
\begin{remark}\label{rem2}\rm
In general, in obvious notation, there are no implications between
$(\r,\s)$ and $(\tilde\r,\tilde\s)$ satisfying Theorem \ref{theo1} (ii), if $\r\ge\tilde\r$. 
Obviously, if Theorem \ref{theo1} (iii) is satisfied for $\r$, it is for $\tilde\r$, if $\r\ge\tilde\r$. 
Note that, if Theorem \ref{theo1} (iii) holds, then $\Phi$ is point-dissipative. 
\end{remark}
We next formulate some well-known definitions, see e.g. \cite{Thieme}. 
\begin{definition}\label{d3}\rm
We call $\ox\in Y$ an equilibrium, if $\Phi(t,\ox)=\ox$ for all $t\ge 0$. Suppose that $E\subset Y$
and that $\ox\in E$ is an equilibrium. We say that $\ox$ is {\it stable on} $E$, if
\[
\forall\;\e>0\;\exists\;\de>0,\;{\rm s.th.}\;d(\Phi(t,x),\ox)<\e\;\forall\;t\ge 0,\;{\rm if}\;d(x,\ox)<\de\;{\rm and}\;x\in E.
\]
We call $\ox$ {\it stable}, if it is stable on $Y$. The equilibrium $\ox$ is called {\it globally attractive}, if $d(\Phi(t,x),\ox)\->0$, 
for all $x\in Y$. If a stable equilibrium is globally attractive it is called
{\it globally asymptotically stable}. 
\end{definition}
Recall that a global attractor need not be stable,
see \cite[Exercise 3.7.7]{Perko} and that a stable equilibrium need not be attractive (consider the zero vector field). 
In \cite[Section 10]{Hartman} for global asymptotic stability of an equilibrium only global existence of
an arbitrary solution along with convergence to the equilibrium is required. 
%
%
\subsection{General differential equations}\label{ss4}
For $n\in\N$ and $h\in(0,\infty)$ we will use the Banach spaces
\begin{eqnarray}
C([-h,0],\R^n),\;\|\phi\|_0:=\max_{\th\in[a,b]}|\phi(\th)|,\;
\nn\\
C^1([-h,0],\R^n),\;\|\phi\|_1:=\|\phi\|_0+\|\phi'\|_0.
\nn
\end{eqnarray}
With $n$ depending on the context, we will refer to the corresponding topologies as $C$- and $C^1$- topologies, respectively, and also denote a function as  $C^1$, if it 
is continuously differentiable. In the context of $N$-dimensional DDE with state-dependent delay (SD-DDE) let us introduce $C:=C([-h,0],\R^N)$ and $C^1:=C^1([-h,0],\R^N)$.  
We will  define solutions as in \cite{Hartung}. Let $U\subset C^1$ be open and $F:U\->\R^N$. 
\begin{definition}\rm
For any $\phi\in U$, a {\it solution} on $[-h,t_*)$ for some $t_*\in(0,\infty]$ of the initial value 
problem (IVP), defined by (\ref{eq55}) and $x_0=\phi$, 
is a $C^1$ function $x:[-h,t_*)\->\R^N$, which satisfies $x_t\in U$ for all 
$t\in(0,t_*)$ as well as the IVP. 
\end{definition}
Solutions on closed intervals $[-h,t_*]$, $t_*>0$, are defined analogously. Admissible initial functions will belong to the set
\begin{eqnarray}
X=X(F):=\{\phi\in U:\;\phi'(0)=F(\phi)\}. 
\label{eq13}
\end{eqnarray}
\begin{definition}\rm
Let ${\cal O}\subset C^1([-h,0],\R^n)$ be open, $f:{\cal O}\->\R^n$. We say that $f$ {\it fulfills
(S)}, if $f$ is $C^1$, each derivative $Df(\phi)$, $\phi\in{\cal O}$ extends to a linear map 
$D_ef(\phi):C([-h,0],\R^n)\->\R^n$ and the map ${\cal O}\times C([-h,0],\R^n)\->\R^n;\;(\phi,\chi)\longmapsto (D_ef)(\phi)\chi$ is continuous.
For arbitrary ${\cal O}_+\subset {\cal O}$, we call $f$ {\it strongly Lipschitz, 
uniformly on bounded sets (sLb) on} ${\cal O}_+$, if for all bounded $B\subset {\cal O}_+$ there exists some $L_B\ge 0$
such that 
\[
|f(\phi)-f(\chi)|\le L_B\|\phi-\chi\|_0\;\forall\;\phi,\chi\in B. 
\]
\end{definition}
Unless stated otherwise, we hereby agree that
when discussing properties of functions, we refer to the topology associated to the domain of the function. In
particular, we will use a state space in $C^1$, such that, for functions defined on a subset of $C^1$, if not 
specifying other topologies, we refer to the $C^1$-topology. In this sense, in the previous definition, a priori
Lipschitz refers to arguments in $C^1$, but since we would like the $C$-norm on the right hand side, 
in \cite{Getto} we introduced the term strongly in this context. 
Suppose that $F$ satisfies (S) with $n:=N$ and ${\cal O}:=U$. A conclusion is the following. 
%
\begin{lemma}\label{lem20}
Suppose that $X\neq\emptyset$. Then $X$ is a $C^1$ sub-manifold  of $U$ with co-dimension $N$. For 
each $\phi\in X$ there exists some $t_\phi$ and a non-continuable global solution $x^\phi:[-h,t_\phi)\->\R^N
$ of the IVP. All segments $x^\phi_t$,
$t\in[0,t_\phi)$, belong to $X$ and the map 
\begin{eqnarray}
&&S:\Ome_S:=\{(t,\phi):\;t\in[0,t_\phi),\;\phi\in X\}\-> X;\;S(t,\phi):=x^\phi_t
\label{eq14}
\end{eqnarray}
defines a continuous semiflow. 
\end{lemma}
The lemma is part of \cite[Theorem 3.2.1]{Hartung}. The 
original result also states the associability of a linear variational equation that refers to the derivatives of
functional and semiflow. To later focus on nonnegative initial conditions we suppose for the remainder of the subsection that there is a certain set $U_+\subset U$, such that for $X_+:=X\cap U_+$ and $\Ome_+:=\{(t,\phi):\;t\in[0,t_\phi),\;\phi\in X_+\}$ one has $X_+\neq \emptyset$ and $S(\Ome_+)\subset X_+$.  Let us denote by $T_\phi:=\{x_t^\phi:\;t\in [0,t_\phi)\}$ the orbit of an 
element $\phi\in X$. 
Next, we suppose that $F$ satisfies (sLb) on $U_+$ and that $T_\phi$ is bounded and $\overline T_\phi\subset U$, whenever $\phi\in U_+$.
Then we can add
\begin{lemma}\label{lem23}
For any $\phi\in X_+$ one has $t_\phi=\infty$. 
\end{lemma}
A variant of this result, where (sLb), boundedness of $T_\phi$ and $\overline T_\phi\subset U$ were 
required on all of the open set $U$, is shown as \cite[Theorem 1.7]{Getto}. The proof of the present lemma 
is analogous and we omit it. Requiring that $F$ satisfies (sLb) on $U_+$, only, will allow to drop the 
requirement of some Lipschitz properties for the functions defining the cell SD-DDE in regard to \cite{Getto}. 
Now,  $S$ induces a continuous semiflow via 
\begin{eqnarray}
S_+:\R_+\times X_+\-> X_+;\;S_+:=S.
\label{eq15}
\end{eqnarray} 
In the following  we elaborate some tools to guarantee Theorem \ref{theo1} (iii) for $Y:=X_+$ and 
$\Phi:=S_+$. 
\begin{lemma}\label{lem21} 
Suppose that 
\begin{eqnarray}
\exists\;K',\;{\it s.th.\;}\forall\;\phi\in X_+\;\exists\;t_0=t_0(\phi)\ge 0,{\it\;s.th.}\;
\|x^{\phi}_t\|_1\le K'\;\forall\; t\ge t_0. 
\label{eq67}
\end{eqnarray}
Then 
\begin{eqnarray}
B:=\{x^\phi_t:\;\phi\in X_+,\;t\ge t_0(\phi)+h\}
\label{eq64}
\end{eqnarray}
 is nonempty and compact. 
\end{lemma}
\Proof 
First, $B\neq\emptyset$ follows from $X_+\neq\emptyset$. 
Define  $\tilde B:=\{x_t^\phi:\;\phi\in X_+,\;t\ge t_0(\phi)\}$. Then $B$ and $\tilde B$ are bounded by $K'$. Let $t\ge t_0+h$, $\th\in[-h,0]$ and $-h\le s\le r\le 0$. Since $\tilde B$
is bounded by $K'$, for any $\phi\in X_+$ the functions $x^\phi$ and $(x^\phi)'$ are bounded by $K'$ on 
$[t_0(\phi)-h,\infty)$. Hence
\[
\sup_{\a\in[t+s+\th,t+r+\th]}\max\{|x^\phi(\a)|,|(x^\phi)'(\a)|\}\le K'\;\forall\;\phi\in X_+.
\]
If one combines this estimate with the mean value theorem one has for some $\xi\in[t+s+\th,t+r+\th]$
that
\begin{eqnarray}
|x_{t+r}^{\phi}(\th)-x_{t+s}^{\phi}(\th)|&=&|x^{\phi}(t+r+\th)-x^{\phi}(t+s+\th)|=
|(x^{\phi})'(\xi)||r-s|
\nn\\
&\le& K|r-s|\;\;\forall\;\phi\in X_+.
\label{eq38}
\end{eqnarray}
For $\th=0$, the estimate yields equicontinuity of $B$. Next, 
\begin{eqnarray}
&&|(x_t^{\phi})'(r)-(x_t^{\phi})'(s)|=|(x^{\phi})'(t+r)-(x^{\phi})'(t+s)|
\nn\\
&=&|F(x^\phi_{t+r})-F(x^\phi_{t+s})|
\le L_{\tilde B}\|x^\phi_{t+r}-x^\phi_{t+s}\|_0
\nn
\end{eqnarray}
for all $\phi\in X_+$ for some $L_{\tilde B}$, since $F$ satisfies (sLb) and $x^\phi_{t+s}$ and $x^\phi_{t+r}$
are elements of the bounded set $\tilde B$. A combination of this estimate with (\ref{eq38})
implies equicontinuity of $B':=\{\phi':\;\phi\in B\}$. 
Equicontinuity of $B$   and $B'$ together with boundedness of $B$ imply compactness of $B$, see e.g. \cite[Lemma 2.1]{Getto}.
\qed

In consistency with the term {\it ultimate boundedness} in Definition \ref{d1}, property 
(\ref{eq67}) is equivalent to point dissipativity of $S_+$. Hence the following lemma could be interpreted as a sufficient criterion, for point dissipativity to be transferable from $C$ to $C^1$. 
\begin{lemma}\label{lem22}
Suppose that $F$ maps $C$-bounded sets into bounded sets and
\begin{eqnarray} 
&&\exists\;K,\;{\it s.th.\;}\forall\;\phi\in X_+\;\exists\;T=T(\phi)\ge 0,{\it\;s.th.}\;
\|x^{\phi}_t\|_0\le K\;\forall\;t\ge T.
\nn\\
\label{eq41}
\end{eqnarray}
Then (\ref{eq67}) holds, i.e., $S_+$ is point dissipative.
%
%
\end{lemma}
\Proof Choose $K$ and $T(\phi)$ according to the assumption. Define $B':=\{\phi\in X_+:\;\|\phi\|_0\le 
K\}$. Then $B'$ is $C$-bounded. Choose $L_B$ such that $|F(\phi)|\le L_B$ for all $\phi\in B'$. Define 
$t_0=t_0(\phi):=T(\phi)+h$ for any $\phi\in X_+$. Then $x^\phi_{t+\th}\in B$ by (\ref{eq41}) for all $\phi\in X_+$, $t\ge t_0$, $\th\in[-h,0]$. Hence, for these values
\begin{eqnarray}
|(x^{\phi}_t)'(\th)|=|(x^{\phi})'(t+\th)|=|F(x^\phi_{t+\th})|\le L_B.
\nn
\end{eqnarray}
Thus,  $\|(x_t^\phi)'\|_0\le L_B$, hence $\|x_t^\phi\|_1\le L_B +K=:K'$, for all $\phi\in X_+$, $t\ge t_0$.
Hence (\ref{eq67}) holds.
\qed

Before we combine the results, let us discuss the conclusion of Lemma \ref{lem21}, e.g., provided
(\ref{eq67}) holds. Obviously, $B$ attracts $X_+$ and Theorem  \ref{theo1} (iii)  is satisfied, irrespective of the choice of $\r$.
If we combine Lemmas \ref{lem21} and \ref{lem22} we get
\begin{corol}\label{corol1}
Suppose that, either $F$ maps $C$-bounded sets into bounded sets and (\ref{eq41}) holds, or (\ref{eq67})
holds. Then Theorem  \ref{theo1} (iii) holds for the set $B$ as in (\ref{eq64}) for all $\phi\in X_+$, irrespective of the choice of $\r$, and $S_+$ is point dissipative. 
\end{corol}
In the remainder of the subsection, we adapt the stability notions of Definition \ref{d3} to the introduced topologies. First, for subsets of $C^1$ one can define the same notions with respect to the metric induced by $\|\cdot\|_0$. Note that a priori there are no implications between (global asymptotic) stability or attractivity with respect to $\|\cdot\|_1$ and the according property with respect to $\|\cdot\|_0$. The following result however shows that under some conditions there are. 
\begin{lemma} \label{lem31} 
Suppose that $0\in X_+$, $F(0)=0$ and that $F|_{U_+}$ is $C$-continuous in zero. 
Then, if for $\phi\in X_+$ and a solution $x=x^\phi$ one has $x_t\->0$ in $C$, then also $x_t\->0$ in $C^1$. In particular, if zero is 
globally attractive in $C$, it is so in $C^1$. Moreover, if for $E\subset X_+$, zero is stable on $E$ in $C$, then also in $C^1$. 
Hence, if zero is globally asymptotically stable in $C$, then also in $C^1$. 
\end{lemma}
\Proof Let $\e>0$. By the continuity assumption for $F$ we can choose $\de$, such that $|F(\phi)|\le\e/2$
if $\|\phi\|_0\le\de$. Next, choose $T$ such that $\|x_t\|_0\le\min\{\de,\e/2\}$ if $t\ge T$. Then for $t\ge T$
one has $|x'(t)|=|F(x_t)|\le\e/2$. Thus for $t\ge T+h$ one has $\|x_t'\|_0\le\e/2$. Hence, $\|x_t\|_1\le\e$ for $t\ge T+h$. We have shown the first statement. The statement on stability can be proven with similar arguments. We omit the details. 
\qed

Recall that (sLb) and (S) require the $C^1$-topology in the domain of the functional. Hence the two 
together do not 
imply $C$-continuity of $F$.
\begin{remark}\rm
Lemma \ref{lem31} can be rewritten for arbitrary constant solutions and the proof is analogous. 
\end{remark}
%
%
%
%
\subsection{The cell FDE}\label{ss2}
We introduce $I:=(R_-,\infty)$, for some $R_-<0$, and consider a function $q:I\->\R$. In the notation of the previous 
subsection, we set $N:=2$ and define $U:=C^1([-h,0],\R)\times C^1([-h,0],I)$. Then, we consider a
functional 
\begin{eqnarray}
j:U\->\R
\label{eq69}
\end{eqnarray}
and a parameter $\mu\ge 0$.
In the following we will analyze the cell FDE. For multiple use, we establish the integrated equations of (\ref{eq68}, \ref{eq58}, \ref{eq70}).
Applying the variation of constants formula to (\ref{eq58}), it becomes clear that solutions, if they exist
correspondingly, satisfy
\begin{eqnarray}
w(t)=e^{\int_0^t q(v(s))ds}\vi(0),
\label{eq12}
\end{eqnarray}
\[
v(t)=\psi(0)e^{-\mu t}+\int_0^t e^{-\mu(t-s)}j(w_s,v_s)ds.
\label{voc} \tag{VOC}
\]
We define 
\begin{eqnarray}
F:U\->\R^2;\;\;F(\vi,\psi):=(q(\psi(0))\vi(0),j(\vi,\psi)-\mu\psi(0)). 
\label{eq63}
\end{eqnarray}
Then, (\ref{eq68}, \ref{eq58}, \ref{eq70}) is, for $x=(w,v)$, of the form (IVP) and $X$ is 
well-defined via (\ref{eq13}, \ref{eq63}). For $\vi:[-h,0]\->\R$, let us write $\vi\ge 0$, if $\vi(\th)\ge 0$ for all $\th\in[-h,0]$. We define 
\begin{eqnarray}
U_+&:=&\{(\vi,\psi)\in U:\;\vi\ge 0,\;\psi\ge 0\}=C^1([-h,0],\R^2_+),\;{\rm which\;yields}
\nn\\
X_+&=&\{\phi\in C^1([-h,0],\R^2_+):\;\phi'(0)=F(\phi)\}.
\label{eq71}
\end{eqnarray}
The following hypotheses suffice to guarantee the existence of a semiflow on $X_+$. 
\begin{hypo}\label{hypo13}
The functional $j$ satisfies (S) on $U$ and (sLb) on $U_+$. Moreover, $j(\vi,\psi)\ge0$
if $(\vi,\psi)\in U_+$, 
and $j(B_1\times B_2)$ is bounded,  whenever $B_1\times B_2\subset U_+$ and $B_1$ is bounded. 
\end{hypo}
\begin{hypo}\label{hypo14}
The function $q$ is bounded and $C^1$. 
\end{hypo}
\begin{hypo} \label{hypo1}
$X_+\neq\emptyset$.
\end{hypo} 
Hence also $X\neq\emptyset$. 
 \begin{remark}\label{rem3}\rm
The existence of everywhere positive elements in $X$ could possibly be proven if one would assume $X\neq\emptyset$ instead of $X_+\neq\emptyset$ and use the specific definition of $X$ and the assumed smoothness properties of $F$. Note that an arbitrary sub-manifold of $C^1([-h,0],\R)$ of co-dimension one need not contain functions that are positive everywhere, as the example $\{\phi:\;\phi(0)=0\}$ shows.
 On the other hand, the functional $F(\phi)=\phi'(0)+\phi(0)$ does not fulfill (S) (it is $C^1$ but the derivative does not allow an extension to $C([-h,0],\R)$), such that the above counter-example does not ``translate" to manifolds generated by functionals via DDE in the discussed way. 
\end{remark}
\begin{lemma}\label{lem39}
The set $X$ is a $C^1$-submanifold of $U$ with co-dimension 
$2$. For each $\phi\in X$ there exists some $t_\phi$ and a non-continuable solution $x^\phi:[0,t_\phi)\->\R^2$ of (\ref{eq68},
\ref{eq58}, \ref{eq70}), all segments $x_t^\phi$ belong to $X$ and the solutions define a continuous 
semiflow $S$ on $\Ome_S\subset\R_+\times X$ via (\ref{eq14}). Moreover, $t_\phi=\infty$ for all $\phi\in X_+$, $S(\Ome_+)\subset X_+$ and thus solutions also define  a continuous semiflow $S_+$ on $\R_+\times X_
+$ via (\ref{eq15}). 
\end{lemma}
\Proof 
For $F=(F_1,F_2)$ (in obvious notation) to see that $F_1$ and $F_2$ satisfy (S), 
one can use evaluation operators, product and composition rules and argue as in the discussion
and proof of \cite[Theorem 1.8]{Getto}. We omit further details on this part of the proof. We can conclude from Lemma \ref{lem20}, that $X$ has the stated 
properties, the existence of non-continuable solutions on $X$ and that a continuous semiflow $S$ is 
defined on $\Ome$ in the discussed way. To guarantee (sLb) for $F$, let $B_1\times B_2\subset U_+$ be 
bounded.  Then 
\[
\overline{conv\{\psi(0):\;\psi\in B_2\}} \subset \R_+,
\]
 where $conv\;A$ denotes the convex hull of a set $A$, is compact
since $B_2$ is $C$-bounded. As 
$q$ is $C^1$, the previous observation and the mean value theorem imply that the functional $\psi\mapsto 
q(\psi(0))$ is (sLb) on $\{\psi:\;(\vi,\psi)\in U_+\}$. Then, it should be 
clear that $F_1$  is (sLb) on $U_+$ and similarly that $F_2$ inherits (sLb) from $j$. 
If one uses the assumed non-negativity property for $j$ and the integrated equations it becomes clear that 
nonnegative initial conditions lead to nonnegative solutions and that $S(\Ome_+)\subset X_+$. This inclusion, as $C^1([-h,0],\R_+^2)$ is (even $C$-)closed, leads to ${\overline T}_\phi\subset C^1([-h,0],\R_+^2)\subset U$ for all $\phi\in X_+$, 
Then, to see the remaining statements it suffices to guarantee the boundedness property for the orbit, as required in the previous subsection, and apply Lemma \ref{lem23}. 
By (\ref{eq12}) and boundedness of $q$, the $w$-component of a finite time orbit is $C$-bounded. Using 
this, (\ref{eq68}) and again boundedness of $q$ one can show directly that the $w$-component of a finite
time orbit is also bounded. With boundedness of the $w$-orbit, the assumed boundedness property for $j$
and (VOC), one can show that the $v$-component of a finite time orbit is $C$-bounded. With previous
arguments and (\ref{eq58}) one can show also boundedness of the $v$-orbit.
\qed

Regarding the establishment of point dissipativity via Lemma \ref{lem22}, under Hypothesis \ref{hypo14}, 
if $j$ maps $C$-bounded sets into bounded sets, then so does $F$. Note that there are no implications 
between $j$ mapping $C$-bounded sets into bounded sets and the boundedness property that was 
assumed for $j$ in Hypothesis \ref{hypo13}. Under the set of assumptions that we will give for the 
ingredients of (\ref{eq68}--\ref{eq5}), however, we will see below that both properties hold. Under 
Hypotheses \ref{hypo13} and \ref{hypo14} we establish the following result on transferability of dissipativity 
for (\ref{eq68}, \ref{eq58}), leaving Lemma \ref{lem22} as a criterion for general DDE 
with $F$ fulfilling the corresponding boundedness property or for (\ref{eq68}--\ref{eq5}) directly. 
\begin{lemma}\label{lem24}
Suppose that there exists some $K>0$, such that for all $\phi\in X_+$ there exists some $T=T(\phi)$, 
such that $w^\phi(t)\le K$ for all $t\ge T(\phi)$. Then for any given persistence function Theorem \ref{theo1} (iii) 
holds and  (\ref{eq68}, \ref{eq58}) is point dissipative in both norms. 
\end{lemma}
\Proof For $\phi\in X_+$ consider $(w,v)=(w,v)^\phi$. Then, $|w'(t)|=|q(v(t))|w(t)\le K\sup|q|$, for all $t\ge T(\phi)$. 
Thus, $\|w_t\|_1\le r:=K(1+\sup|q|)$, for all $t\ge T(\phi)+h$. Then, in obvious notation, by (VOC)
for all $t\ge T(\phi)+h$
\begin{eqnarray}
v(t)&=&e^{-\mu(t-(T(\phi)+h))}v(T(\phi)+h)+e^{-\mu t}\int_{T(\phi)+h}^tj(w_s,v_s)ds
\nn\\
&=:&(I)+(II).
\nn
\end{eqnarray}
We can choose $t_1=t_1(\phi)\ge T(\phi)+h$ so large that $(I)\le 1$ for all $t\ge t_1$. Next, by Hypothesis
\ref{hypo13} we can choose $L_B$, such that $j(\vi,\psi)\le L_B$ whenever $\vi\in \overline B_r(0)$ and 
$(\vi,\psi)\in U_+$. 
Then $(II)\le L_B/\mu$ for all $t\ge t_1$. Thus, for $K_1:=\max\{K,1+L_B/\mu\}$ and $t_1(\phi)$
ultimate boundedness in the $C$-norm follows. Finally, 
\begin{eqnarray}
|v'(t)|\le j(w_t,v_t)+\mu v(t)\le L_B+\mu(1+\frac{L_B}{\mu})=2L_B+\mu\;\forall\;t\ge t_1(\phi).
\nn
\end{eqnarray}
This implies ultimate boundedness in the $C^1$-norm, i.e., that (\ref{eq67}) holds. The remaining statements follow by Corollary \ref{corol1}. 
\qed

Part (b) of the following result (transferability of stability from the $w$-component to the system) will be used
in our proof of stability of zero below. To establish convergence to zero we will consider
the case, where $q(0)<0$ is given, separately, since it saves some assumptions that are needed otherwise. Then, part (c) (transferability of convergence to zero from the $w$-component to the system) will be used to establish convergence to zero in case $q(0)<0$. For the case where $q(0)<0$ is not given, it will turn out convenient to first establish that $v_t\->0$ and we will see that part (c) is not necessary. It will be used in the persistence proof.
\begin{lemma}\label{lem8} 
\begin{itemize}
\item[(a)] Equivalently $j(0)=0$,  $F(0)=0$, $0\in X_+$ or (\ref{eq68}, \ref{eq58}, \ref{eq70}) has a zero 
solution. In either case Hypothesis \ref{hypo1} holds. If $j|_{U_+}$ is $C$-continuous in zero, then $F|_{U_+}$ is $C$-continuous in zero. Hence, if both 
$j(0)=0$ and $j|_{U_+}$ is $C$-continuous in zero, then the conclusions of Lemma \ref{lem31} hold. 
\end{itemize}
\noindent
If both, $j(0)=0$ and the $C$-continuity of $j$ in zero is uniform with respect to the second argument, i.e., if
\begin{eqnarray}
&&\forall\;\e>0\;\exists\;\de>0,\;{\it s.th.}\;j(\vi,\psi)\le\e\;
\forall\;(\vi,\psi)\in U_+\;with\;
\|\vi\|_0\le\de,
\nn\\
\label{eq19}
\end{eqnarray}
then also the following hold. 
\begin{itemize}
\item[(b)] Let $E\subset X_+$ with $0\in E$. Suppose that $\|w^{(\vi,\psi)}_t\|_0\le\|\vi\|_0$ for all $(\vi,\psi)\in E$, $t\ge0$. Then zero is stable on $E$ in $C$ and in $C^1$. 
\item[(c)] Let $\phi\in X_+$ and suppose that for $(w,v)=(w,v)^\phi$ one has $w\->0$  as $t\_>\infty$. Then  $(w,v)_t
\->0$ in $C$ and in $C^1$. In particular, if on $X_+$ one has $w\->0$, then zero is globally attractive in  
both norms. 
\end{itemize}
\end{lemma}
\Proof (a) The statements on equivalence are clear. Next, note that the evaluation operator $\phi\longmapsto\phi(0)$ defined on a subset of $C^1([-h,0],\R^n)$ is $C$-continuous and recall that $q$ is continuous. Then it is easy to see that also the continuity statement for $F$ holds. 
\\
(b) Let $\e>0$. By the continuity assumption for $j$ we can choose $\de\le\e/2$, such that 
$j(\vi,\psi)\le\e\mu/2$ if $(\vi,\psi)\in E$ with $\|\vi\|_0\le \de$. Let  $(\vi,\psi)\in E$ with 
$\|(\vi,\psi)\|_0\le \de$. Then for $(w,v)=(w,v)^{(\vi,\psi)}$ one has $\|w_t\|_0\le\|\vi\|_0\le\de\le\e$. Moreover
by (VOC) the above estimate for $j$ yields for $t\ge0$
\begin{eqnarray}
v(t)&\le& e^{-\mu t}\psi(0)+\frac{\e}{2}(1-e^{-\mu t})\le\psi(0)+\frac{\e}{2}\le\de+\frac{\e}{2}\le\e. 
\nn
\end{eqnarray}
This implies the stated $C$-stability and $C^1$-stability follows from (a). 
\\
(c) Choose $(\vi,\psi)\in X_+$ and let $\e>0$. Choose $\de\le\e/2$ such that $j(\vi,\psi)\le\e\mu/2$
if $\|\vi\|_0\le\de$. Choose $t_1$ such that $\|w_t\|_0\le\de$ if $t\ge t_1$. Then $\|w_t\|_0\le\e$ for all
$t\ge t_1$. Choose $t_2\ge t_1$ such that 
\[
e^{-\mu t}[\psi(0)+\int_0^{t_1}e^{\mu s}j(w_s,v_s)ds]\le\e/2\;\forall\;t\ge t_2.
\]
Then for $t\ge t_2$ by (VOC)
\begin{eqnarray}
v(t)&=&e^{-\mu t}[\psi(0)+\int_0^{t_1}e^{\mu s}j(w_s,v_s)ds]+
e^{-\mu t}\int_{t_1}^te^{\mu s}j(w_s,v_s)ds
\nn\\
&\le&\e/2+\e/2=\e. 
\nn
\end{eqnarray}
 Hence, for $t\ge t_2+h$, one has $\|(w,v)_t\|_0\le\e$. 
We have shown that $\|(w,v)_t\|_0\->0$, as $t\_>\infty$. Again, the statement on the $C^1$-norm follows by (a). 
\qed

%
\section{GAS of zero for the cell FDE}\label{s2}
We consider (\ref{eq68}, \ref{eq58}, \ref{eq70}) and keep Hypotheses \ref{hypo13} -- 
\ref{hypo1}
and the notation of Subsection \ref{ss2} throughout the section. 
We show that under certain additional conditions
there is a globally asymptotically stable trivial equilibrium. 
\begin{lemma}\label{lem9}
Suppose that $q(z)\le 0$, for all $z\in[0,\infty)$ and $\phi\in X_+$. Then for $w=w^\phi$ one has  
$w'(t)\le 0$ for all $t\ge 0$ and there
exists some $w^\infty=w^\infty(\phi)\in[0,\infty)$ such that $w(t)\->w^\infty$ as $t\_>\infty$. 
If moreover (\ref{eq19}) holds, then by Lemma \ref{lem8}  there is a zero solution which is stable 
on $X_+$ in $C$ and in $C^1$. 
\end{lemma}
\Proof Since $w(t)\ge 0$ and $v(t)\ge 0$ for all $t\ge 0$ we have by (\ref{eq68}) that
$w'(t)=q(v(t))w(t)\le 0$ for all $t\ge 0$ and the first statement is proven. The convergence statement
follows since $w$ is non-increasing and nonnegative. The statement on stability follows by Lemma \ref{lem8}  (b) applied to $E=X_+$. 
\qed

\begin{hypo}\label{hypo15}%
One has $q(s)<0$ for all $s>0$. 
\end{hypo}
The hypothesis allows for $q(0)=0$ and for the specification of $q$ with the 
parameters given in Section \ref{s6} below. 
There cannot be positive zeros of $q$, hence there cannot be a positive equilibrium. 

Let $(w,v)$ be a solution. Then $w'\le 0$. Moreover, if $\psi(0)>0$, then $v(s)>0$ for all $s\ge0$ and hence,
if also $\vi(0)>0$, then $w'<0$. 
and $w(t)\->w_\infty$ as $t\_>\infty$ for some $w_\infty\ge0$. In particular, $w$ is bounded on $\R_+$. Since $q$ is bounded, $w'$ is bounded on $\R_+$. 
By Hypothesis \ref{hypo13}, the set $\{j(w_t,v_t):\;t\in\R_+\}\subset\R$ is also bounded. It follows from (VOC)
that $v$ is bounded on $\R_+$. %
\begin{lemma}\label{lem25}
The functions $w$, $w'$, $v$ and $v'$ are bounded on $\R_+$, $w(t)\->w_\infty$ as $t\_>\infty$ for some 
$w_\infty\ge0$ and $w'(t)\->0$ as $t\_>\infty$. 
\end{lemma}
\Proof All assertions have been shown except for the last one. According to \cite[Corollary A 17]{Thieme}
it is sufficient to show that $|w'(t+s)-w'(t)|\->0$ as $t\_>\infty$ and $s\_>0$. By (\ref{eq68})
\begin{eqnarray}
&&|w'(t+s)-w'(t)|
\nn\\
&\le&|q(v(t+s))-q(v(t))||w(t+s)|+|q(v(t))||w(t+s)-w(t)|. 
\nn
\end{eqnarray}
Since $q\circ v$ is bounded and $w(t)\->w_\infty$, the second term on the left hand side converges to zero
as $t\_>\infty$ and $s\_>0$. Since $w$ is bounded, it is sufficient to show that 
\[
|q(v(t+s))-q(v(t))|\->0,\;\;{\rm as}\;t\_>\infty,\;s\_>0. 
\]
Recall that $v$ is bounded. Since $q$ is uniformly continuous on any bounded subset of $\R_+$, it is sufficient to show that $|v(t+s)-v(t)|\->0$ as $t\_>\infty$ and $s\_>0$. By the mean value theorem, for some $r\in[t,t+s]$
one has 
\[
|v(t+s)-v(t)|=|v'(r)|s. 
\]
Since $v'$ is bounded, this tends to zero, as $s\_>0$. 
\qed

\begin{lemma}\label{lem26}
One has $w(t)\->0$ or $v(t)\->0$ as $t\_>\infty$. 
\end{lemma}
\Proof Assume that $w_\infty>0$. Since $w'(t)\->0$ by Lemma \ref{lem25}, one has $q(v(t))\->0$ by (\ref{eq68}). 
Suppose that $v$ does not converge to zero. Since $v$ is bounded, there exists a sequence $(t_n)$ and some $\tilde v>0$, such that $t_n\_>\infty$ and $v(t_n)\->\tilde v$ as $n\_>\infty$. Then $q(v(t_n))\-> q(\tilde v)<0$ as $n\_>\infty$, a contradiction. 
\qed

Suppose that for the remainder of the section (\ref{eq19}) holds. 
\begin{prop}\label{p1}
One has $v(t)\->0$ and $v'(t)\->0$ as $t\_>\infty$.
\end{prop}
\Proof By Lemma \ref{lem26}, we can assume that $w\->0$. Since $v$ is bounded, by \cite[Proposition A 22]{Thieme}
there exists a sequence $(t_n)$ with $t_n\_>\infty$ and 
\[
v(t_n)\->v^\infty:=\limsup_{t\_>\infty}v(t)
\]
and $v'(t_n)\->0$ as $n\_>\infty$. By (\ref{eq58})
\[
0=\lim_{n\_>\infty}j(w_{t_n},v_{t_n})-\mu v^\infty. 
\]
Since $\|w_{t_n}\|\->0$ as $n\_>\infty$, the first term on the right hand side equals zero by (\ref{eq19}). This implies that $v^\infty=0$ and $v\->0$. By (\ref{eq19}) and (\ref{eq58}), one concludes that $v'\->0$. 
\qed

In the following we distinguish the case where $q(0)<0$ from the case, where this is not given, starting 
with the former. 
\begin{lemma}\label{lem27}
If $q(0)<0$, then $w(t)\->0$ as $t\_>\infty$. 
\end{lemma}
\Proof One has $q(0)\le -\de<0$ for some $\de>0$. Choose $z_1>0$ such that $q(z)\le-\de/2$ for all 
$z\in[0,z_1]$. By Proposition \ref{p1} one has $z\->0$ and we can choose $T$ such that $z(t)\in[0,z_1]$
for all $t\in[T,\infty)$. Then for some $K>0$ and $t\ge T$
\begin{eqnarray}
w(t)=e^{\int_0^tq(v(s))ds}\vi(0)\le K e^{\int_T^tq(v(s))ds}\vi(0)\le K\vi(0)e^{-\frac{\de}{2}(t-T)}\->0.
\nn
\end{eqnarray}
\qed

Now, if $q(0)<0$, GAS of zero can be concluded as follows. First recall that since (\ref{eq19}) holds, stability of zero follows by Lemma \ref{lem9} (which relies on Lemma \ref{lem8} (b)) and note that global attractivity of zero in both norms can be concluded from the previous result and  Lemma \ref{lem8} (c). Without quoting Lemma \ref{lem8}, attractivity of zero in both norms can be concluded from Lemma \ref{lem27}, boundedness of $q$ and (\ref{eq68}) ($w$-component) and Proposition \ref{p1} ($v$-component).   
We thus get:
\begin{corol}
If $q(0)<0$, then zero is GAS. 
\end{corol}
For the case, where $q(0)<0$ is not given, we formulate
\begin{hypo}\label{hypo9}
One has $j(\vi,0)>0$ if $\vi$ is a strictly positive constant function. 
\end{hypo}
\begin{prop}
One has $x_t\->0$ as $t\_>\infty$ in both norms. 
\end{prop}
\Proof By our previous results, we have that $w\->w_\infty$, $w'\->0$, $v\->0$ and $v'\->0$. It remains to show that
$w_\infty=0$. By (\ref{eq58}), one has $j(w_t,v_t)\->0$. Note that, in the $C^1$-norm, $(w,v)_t\->(w_\infty,0)$. 
Then, by continuity of $j$, one has that $j(w_t,v_t)\->j(w_\infty,0)$. Hence, if $w_\infty>0$, by Hypothesis 
\ref{hypo9} we get $j(w_\infty,0)>0$, which contradicts $j(w_t,v_t)\->0$. 
\qed

We can summarise the section as follows. 
\begin{theorem}\label{theo5}
Suppose that Hypotheses \ref{hypo13} -- \ref{hypo14} as well as (\ref{eq19}) hold and consider 
(\ref{eq68}, \ref{eq58}, \ref{eq70}). Then zero is stable. If moreover Hypothesis \ref{hypo15} and, at least one of the two, $q(0)<0$ or Hypothesis \ref{hypo9}, hold, then zero is GAS. 
\end{theorem}
\begin{remark}\rm
Recall that the differentiability assumptions on $q$ and $j$ are used to establish well-posedness. We remark without proof that, under the assumption of well-posedness, continuity properties of $q$ and $j$ would
be sufficient for the results of this section. 
\end{remark}
%
%
\section{Persistence and point-dissipativity for the cell FDE}
\label{s4}
Throughout the section, we keep the setting of Subsection \ref{ss2} and suppose that Hypotheses 
\ref{hypo13} -- \ref{hypo1} hold. 
 We will apply results of Subsection \ref{ss3} and will use some notation introduced there. 
 %
%
%
\subsection{Uniform weak persistence for $\r(\vi,\psi)=\vi(0)$}\label{s3}
In the following, we would like to define persistence functions. By definition, a persistence function is non-identical zero. Therefore we would like to guarantee that $X_+$ contains certain non-trivial elements. 
\begin{hypo}\label{hypo8}
There exists some $(\vi,\psi)\in X_+$ for which $\vi(0)>0$. 
\end{hypo} 
The hypothesis obviously sharpens Hypothesis \ref{hypo1}. 
Then, in relation to Subsection \ref{ss3}, we set $Y:=X_+$ and $\Phi:=S_+$ and define a continuous functional on $X_+$ via (\ref{eq16}). 
%
%
By Hypothesis \ref{hypo8} the functional is non-identical zero. Hence, $\r_1$ is a persistence function.  
\begin{hypo}\label{hypo16}
One has $q(0)>0$.
\end{hypo}
\begin{theorem}\label{theo4}
(uniform weak $\r_1$-persistence) Suppose that Hypotheses \ref{hypo13} -- \ref{hypo14}, (\ref{eq19}) and 
Hypotheses \ref{hypo8} and \ref{hypo16} hold. Define $\r_1$ as in (\ref{eq16}), $F$ as in (\ref{eq63}) and 
$X_+$ as in (\ref{eq71}). Then the semiflow associated to (\ref{eq68}, \ref{eq58}) on $X_+$ in 
Lemma \ref{lem39} is uniformly weakly $\r_1$-persistent. 
\end{theorem}

\Proof First recall that Lemma \ref{lem39} relies on Hypotheses \ref{hypo13} -- \ref{hypo1} and  that
Hypothesis \ref{hypo1} and existence of $\r_1$ are implied by Hypothesis \ref{hypo8}. Next, by Hypothesis \ref{hypo16} we can choose
$\e_0>0$ and $q_1>0$, such that $q(s)\ge q_1$ for all $s\in[0,\e_0]$. Define $\e:=\e_0\mu/2$ and choose 
$\de>0$ according to (\ref{eq19}). Suppose that $S_+$ is not uniformly weakly $\r_1$-persistent. Then there 
exists some $(\vi,\psi)\in X_+$ with $\vi(0)>0$ and a solution $(w,v)=(w,v)^{(\vi,\psi)}$, such that
\[
\limsup_{t\_>\infty}w(t)<\de.
\]
Since $\vi(0)=w(0)>0$, we have $w(t)>0$ for all $t\ge 0$. 
Then we can choose some $T>0$, such that 
$0<w(t)\le \de$ for all $t\ge T$. Hence, $\|w_t\|_0\le \de$ for all $t\ge T+h$. Thus, by (\ref{eq19}) and (VOC)
for all $t\ge T+h$ we get
\begin{eqnarray}
&&v(t)=v(T+h)e^{-\mu [t-(T+h)]}+\int_{T+h}^tj(w_s,v_s)e^{-\mu(t-s)}ds\;{\rm and}
\nn\\
&&
\int_{T+h}^tj(w_s,v_s)e^{-\mu(t-s)}ds\le\e\int_{T+h}^te^{-\mu(t-s)}ds
=\frac{\e}{\mu}(1-e^{-\mu [t-(T+h)]})
\nn\\&&\le\frac{\e}{\mu}=\frac{\e_0}{2}.
\nn
\end{eqnarray}
It follows that $v(t)<\e_0$ for all $t\ge t_0$ and some $t_0>0$. Thus
\[
w(t)=w(t_0)e^{\int_{t_0}^tq(v(s))ds}\ge w(t_0)e^{q_1(t-t_0)}\->\infty,
\]
as $t\_>\infty$. This is a contradiction. 
\qed

%
%
\subsection{Point dissipativity}
We first add the following hypotheses to our assumptions. 
\begin{hypo}\label{hypo17}
One has
\begin{eqnarray}
\exists\;\de>0,\;z^*\ge0, \;{\it s.th.}\;q(z)\le-\de \;\forall\;z\ge z^*.
\nn
\end{eqnarray}
\end{hypo}
\begin{hypo}\label{hypo3}
\begin{eqnarray}
\forall\;K>0\;\exists\;L=L(K),\;{\it s.th.}\;j(\vi,\psi)\ge K,
{\it if}\;(\vi,\psi)\in U_+,\;\min\vi\ge L.
\nn
\end{eqnarray}
\end{hypo}
\begin{lemma}\label{lem13}
There exist $K_1,\e_w>0$ and $t_d$ (``d" for ``decrease") such that the following holds.  
Suppose that for $\phi\in X_+$ there exist $T_0=T_0(\phi)$ and $T=T(\phi)>t_d+h$ such that for $w=w^\phi
$ one has $w\ge K_1$ on $[T_0,T_0+T]$, then $w'\le-\e_w$ on $[T_0+h+ t_d,T_0+T]$. 
\end{lemma}
\Proof First note that by (VOC) for $t\ge T_0+h$
\begin{eqnarray}
v(t)&\ge&\int_{T_0+h}^te^{-\mu(t-s)}j(w_s,v_s)ds\ge\frac{1-e^{-\mu[t-(T_0+h)]}}{\mu}\inf_{s\in[T_0+h,t]} j(w_s,v_s).
\nn
\end{eqnarray}
Now, define $K:=(y^*+1)\mu$ and choose $L=L(K)$ according to Hypothesis \ref{hypo3}. Define $K_1:=L$. 
Let $t\in[T_0+h, T_0+T]$, $s\in[T_0+h,t]$. Then by assumption $w_s(\th)=w(s+\th)\ge K_1=L$. Thus by Hypothesis \ref{hypo3}
\[
\frac{1-e^{-\mu[t-(T_0+h)]}}{\mu}\inf_{s\in[T_0+h,t]} j(w_s,v_s)\ge(1-e^{-\mu[t-(T_0+h)]})(y^*+1).
\]
Hence $v(t)\ge y^*$ if $t\in[T_0+h+t_d,T_0+T]$ if we define 
$t_d=\frac{1}{\mu}\ln(y^*+1)$. 
Hence, $q(v(t))\le-\de$ for all $t\in[T_0+h+t_d,T_0+T]$. Thus for these values
$w'(t)=q(v(t))w(t)\le-\de w(t)\le-\de K_1\le-\e_w$, if we choose $\e_w\le K_1\de$.
\qed

\begin{lemma}\label{lem14}
Choose $K_1$ according to Lemma  \ref{lem13}. 
Then there exists some $K_2>K_1$, such that for all $\phi\in X_+$ the following holds.  Suppose that there exists some  $t_0=t_0(\phi)$ such that  for $w=w^\phi$ one has that $w(t_0)\le K_1$. Then $w(t)<K_2$ for all $t\ge t_0$. 
\end{lemma}
\Proof 
Choose also $t_d$ and $\e_w$ according to Lemma \ref{lem13}. 
Choose $K_2$, such that 
\begin{eqnarray}
K_2> K_1e^{q_0(t_d+h)},\;q_0:=\max\{1,\sup_{z\in[0,z^*]}q(z)\}. 
\label{eq45}
\end{eqnarray}
Suppose that the statement is not true. Then there exist $T_0$, $t_1$, $T_0<t_1$, such 
that $w(T_0)=K_1$, 
$w(t_1)=K_2$ and $w(t)\in[K_1,K_2)$ for all $t\in[T_0,t_1)$. Then
\begin{eqnarray}
K_2=w(t_1)=w(T_0)e^{\int_{T_0}^{t_1}q(v(s))ds}\le K_1e^{q_0(t_1-T_0)}.
\nonumber
\end{eqnarray}
This and (\ref{eq45}) imply that
\begin{eqnarray}
t_1-T_0\ge\frac{1}{q_0}\ln\frac{K_2}{K_1}>t_d+h.
\nn
\end{eqnarray}
Hence by Lemma \ref{lem13} (applied to $T:=t_1-T_0$) one has $w'\le -\e_w<0$
on $[T_0+h+t_d,t_1]$. But then one deduces the contradiction $K_2=w(t_1)<w(T_0+h+ t_d)<K_2$.
\qed

\begin{lemma}\label{lem16}
There exists $K_2$ such that for any $\phi\in X_+$ there exists some $t_m=t_m(\phi)$ such that for $w=w^\phi$ one has that  $w\le K_2$ on $[t_m,\infty)$. 
\end{lemma}
\Proof 
Choose $K_1$, $t_d$ and $\e_w$ according to Lemma \ref{lem13}, $K_2$ according to 
Lemma \ref{lem14}. Fix $\phi$ and consider $w=w^\phi$. There cannot be more than two
cases:
\\
{\it Case 1:} There exists some $t_0=t_0(\phi)$ such that $w(t_0)<K_1$. Then by Lemma \ref{lem14} the statement follows for $t_m:=t_0$. 
\\
{\it Case 2:} One has $w\ge K_1$ on $[T_0,\infty)$ for some $T_0>0$. Then by Lemma \ref{lem13} for any $T>h+t_d$ one has
$w'\le-\e_w$ on $[T_0+h+t_d, T_0+T]$. Hence $w'\le-\e_w$ 
on  $[T_0+h+t_d, \infty)$. This contradicts nonnegativity of $w$.  
\qed

We can now combine Lemma \ref{lem24} with the previous result and conclude the subsection. 
\begin{theorem}\label{theo6}
(point dissipativity)
Suppose that Hypotheses \ref{hypo13} -- \ref{hypo1} and Hypotheses \ref{hypo17} --
\ref{hypo3} hold. Then the semiflow $S_+$ associated to (\ref{eq68}, \ref{eq58}) in Lemma \ref{lem39} satisfies 
Theorem \ref{theo1} (iii) for any given persistence function and is point dissipative (ultimately bounded) in the $C$ 
and the $C^1$ norm. 
\end{theorem}
\begin{remark} \rm
In general, other dissipativity notions can be defined using persistence functions, such as $\r$-disspativity, 
see \cite[Definition 3.2]{Smith}. If Hypothesis \ref{hypo8} holds, for our choices $\r=\r_1$ and $\r=\r_m$ point dissipativity of the semiflow induced by (\ref{eq68}, \ref{eq58}) clearly implies $\r$-dissipativity.
\end{remark}
%
\subsection{Uniform weak persistence via repellance and acyclicity, uniform strong persistence, point dissipativity
and positive equilibrium}
To define a second persistence function, we sharpen Hypotheses \ref{hypo1}  and \ref{hypo8}:
\begin{hypo}\label{hypo10}
There exists $(\vi,\psi)\in X_+$ for which $\min\{\vi(0),\psi(0)\}>0$. 
\end{hypo}
Then we define a continuous functional on $X_+$ via (\ref{eq10}). 
Note that $\r_1(\phi)\ge\r_m(\phi)$ for all $\phi\in X_+$. 
We assume that Hypothesis \ref{hypo16} holds for the remainder of the subsection. 
\begin{lemma}\label{lem1}
Suppose that $\phi=(\vi,\psi)\in X_+$ with $\vi(0)> 0$. Then, for $(w,v)=(w,v)^\phi$ one has that 
$v\->0$ implies $w\-> \infty$. In particular, $(w,v)$ does not converge to zero at infinity. 
\end{lemma}

\Proof If $v\->0$, then $q(v(t))\->q(0)$. Hence we can choose some $T$, such that $q(v(t))\ge \frac{q(0)}{2}>0$ for all
$t\ge T$. Then the statement follows, as for $t\ge T$
\begin{eqnarray}
w(t)&=&\vi(0)e^{\int_0^Tq(v(s))ds}e^{\int_T^tq(v(s))ds}\ge
\vi(0)e^{\int_0^Tq(v(s))ds}e^{\frac{q(0)}{2}(t-T)}.
\nn
\end{eqnarray}
\qed

\begin{remark}\label{rem1}
$\psi(0)>0$ is not sufficient for the previous result: If $\vi(0)=0$, then $w(t)=0$ for all $t\ge0$. Then, if 
(\ref{eq19}) holds, by Lemma \ref{lem8} (c) one has $(w,v)\->0$. 
\end{remark}
Next we ensure that $\{(0,0)\}\subset X_+$ by making:
\begin{hypo}\label{hypo7}
One has $j(0)=0$. 
\end{hypo}
Then, in relation to Subsection \ref{ss3}, we set  $k:=1$ and $M:=\{(0,0)\}$. 
Lemmas \ref{lem1} and \ref{lem17} immediately imply
\begin{lemma}\label{lem34}
The set $M=\{(0,0)\}\subset X_+$ is weakly $\rho_1$-repelling, hence also weakly $\r_m$-repelling.
\end{lemma}
The following hypothesis excludes a trivial $j$ and suffices to give $X_0$ and $\Omega$ a shape that
allows to apply the criterion for weak persistence. 
\begin{hypo}\label{hypo5}
\begin{eqnarray}
\exists\;K, \;{\it s.th.}\;j(\vi,0)>0,\;{\it if}\;(\vi,\psi)\in U_+,\;\min \vi\ge K. 
\nn
\end{eqnarray}
\end{hypo}
\begin{lemma}\label{lem19}
One has $X_0^{\r_1}=X_0^{\r_m}=\{(\vi,\psi)\in X_+:\;\vi(0)=0\}$. 
\end{lemma}
\Proof Define 
\[
B:=\{(\vi,\psi)\in X_+:\;\vi(0)=0\}.
\]
Since $X_0^{\r_1}\subset X_0^{\r_m}$, see Subsection \ref{ss3}, it is sufficient to show that
$X_0^{\r_m}\subset B\subset X_0^{\r_1}$. To show the first inclusion, let
\begin{eqnarray}
\phi=(\vi,\psi)\in X_0^{\r_m}=\{\phi\in X_+:\;\min\{|w^\phi(t)|,|v^\phi(t)|\}=0\;\forall\;t\ge 0\}. 
\nn
\end{eqnarray}
and consider $(w,v)=(w,v)^\phi$. Since $\min\{|w(t)|,|v(t)|\}=0$ for all $t\ge 0$, we have $\vi(0)=0$ or $\psi(0)=0$
or both. If $\vi(0)=0$ we are done. If $\vi(0)\neq0$, then $\psi(0)=0$ and $\vi(0)>0$.
Then $w(t)>0$ for all $t\ge0$. Since $\min\{|w(t)|,|v(t)|\}=0$ for all $t\ge 0$ it must be that $v(t)=0$ for all
$t\ge 0$. Then by (\ref{voc})
\[
\int_0^te^{\mu s}j(w_s,v_s)ds=0\;\forall\; t\ge 0.
\]
Thus, for all $t\ge 0$ one has $j(w_t,v_t)=0$. Hence, $j(w_t,0)=0$ for all $t\ge h$. 
Choose $K$ according to Hypothesis \ref{hypo5}. Since $w(t)=\vi(0)e^{q(0)t}$ and $q(0)>0$ we can 
choose $t\ge h$ such that $w_t(\th)\ge K$ for all $\th\in[-h,0]$. Then by Hypothesis \ref{hypo5} one has
$j(w_t,0)>0$, which is a contradiction. We have shown that $X_0^{\r_m}\subset B$. Now, clearly
\begin{eqnarray}
X_0^{\r_1}=\{\phi\in X_+:\;w^\phi(t)=0\;\forall \;t\ge 0\}. 
\nn
\end{eqnarray}
Let $\phi=(\vi,\psi)\in B$. Then $\vi(0)=0$. Thus by (\ref{eq12}) $w^\phi(t)=0$ for all
$t\ge 0$, hence $\phi\in X_0^{\r_1}$. This completes the proof. 
\qed

Since $X_0^{\r_1}=X_0^{\r_m}$, from now on we simply write 
$X_0$. Obviously, $M\subset X_0$. 
Moreover, clearly $M$ is isolated in $X_+$, compact and invariant. 
We next guarantee that $\Ome\subset M$. Note that in our setting, the closure in Definition \ref{d2} 
refers to the $C^1$-topology. 
The following result relates to this definition.  Its proof and the proofs of the remaining results of the subsection refer to  Lemma \ref{lem8} (b) and (c). To apply these results, we sharpen Hypothesis \ref{hypo7} by assuming that, for the remainder of the subsection, (\ref{eq19}) holds. 
\begin{lemma}\label{lem38}
Let $x\in X_0$. Then $\ome(x)=\{(0,0)\}$, hence also $\Ome=\{(0,0)\}$ and thus $\Ome=M$.
\end{lemma}
\Proof Since $x\in X_0$, by Lemma \ref{lem8} (c) we have
\begin{eqnarray}
S(t,x)\->0,\;{\rm for}\; t\_>\infty.
\label{eq35}
\end{eqnarray} 
``$\supset$'' We should show that for any $t\ge 0$ one has
$
(0,0)\in\overline{S([t,\infty)\times x)}
=\overline{\{S(s,x):\;s\in[t,\infty)\}}.
$
This follows from (\ref{eq35}). 
\\
``$\subset$'' Let $y\in\ome(x)$. Then for all $t\ge 0$ there exists a sequence
$(y_n^t)\in (S([t,\infty)\times x))^\N$ such that $y_n^t\->y$ as $n\_>\infty$. 
Hence for all $t\ge 0$ there exists a sequence
$(s_n^t)\in [t,\infty)^\N$ such that $S(s_n^t,x)=y_n^t\->y$ as $n\_>\infty$. 
Thus for all $j\in\N$ one has $(s_n^j)\in [j,\infty)^\N$ and $S(s_n^j,x)\->y$ as $n\_>\infty$. 
Then for all $j\in\N$ there exists some $s_{n_j}^j\in[j,\infty)$ with
\[
\|S(s_{n_j}^j,x)-y\|_1\le\frac{1}{j}.
\]
Define $r_j:=s_{n_j}^j$. Then $r_j\->\infty$ and $\|S(r_j,x)-y\|_1\->0$, as $j\_>\infty$.
Now, $y=0$ by (\ref{eq35}). Thus $\ome(x)\subset \{(0,0)\}$ and ``$\subset$"
is shown. 
\qed

Lemma \ref{lem19} and Lemma \ref{lem8} (b) applied to $E=X_0$ immediately imply
\begin{lemma}\label{lem18}
$(0,0)$ is stable on $X_0$, in both, $C$ and $C^1$-norm.
\end{lemma}
In a metric space $Y$, let us denote by $B_\e(x_0)$ a ball with radius $\e>0$ around $x_0\in Y$. 
\begin{lemma}
$\{(0,0)\}$ is acyclic. 
\end{lemma}
\Proof Suppose that $\{(0,0)\}$ is cyclic. Then there exists a total trajectory $T:\R\->X_0$, such that
$T(0)\neq (0,0)$ and $T(-t)\->\{(0,0)\}$ as $t\rightarrow\infty$. Hence, $S_+(t,T(s))=T(t+s)$, for all $t\ge 0$, 
$s\in \R$ and
\begin{eqnarray}
&&\forall\;U\subset X_0\;{\rm open}\;, \{(0,0)\}\subset U\;\exists\;r\in\R_+, \;{\rm s.th.} \;T(-t)\in U\;\forall\;t\in[r,\infty).
\nn\\
\label{eq18}
\end{eqnarray}

As $T(0)\neq (0,0)$, we have $\e:=(1/2)\|T(0)-(0,0)\|_1>0$  and $T(0)\notin B_\e(0,0)$.
By Lemma \ref{lem18} we can choose $\de$ such that $x^\phi_t\in B_\e(0,0)$ for all $t\ge 0$
if $\phi\in B_\de(0,0)\cap X_0$. Next, by (\ref{eq18}) we can choose some $\ot>0$ such that $T(-\ot)\in  B_\de(0,0)\cap X_0$. Then, $T(0)=S_+(\ot,T(-\ot))=x^{T(-\ot)}_{\ot}\in  B_\e(0,0)$, which is a contradiction.
\qed

In the following, we combine the results on repellence and acyclicity of this subsection with the compactness 
results of the previous subsection and the theory of Section \ref{ss3} to conclude uniform weak persistence. 
Uniform persistence and the existence a positive equilibrium follow without further assumptions. 

For $\s_1:=\r_1\circ S_+$ and $\s_m:=\r_m\circ S_+$  one has 
\begin{eqnarray}
\s_1(t,\phi)=w^\phi(t),\;\s_m(t,\phi)=\min\{w^\phi(t),v^\phi(t)\}.
\end{eqnarray}
The functions $\s_1$ and $\s_m$ are continuous as compositions of continuous functions, hence Theorem
\ref{theo1} (i) holds for both. In a similar way as we have shown the preservation 
of non-negativity in Subsection \ref{ss2} it can now be shown that $S_+$ satisfies also Theorem \ref{theo1} 
(ii) for all $y\in X_+$, i.e., in particular irrespective of the choice of $B$, and with respect to both persistence functions. We omit the details. 
\begin{theorem}\label{theo2} (uniform persistence, point dissipativity and positive equilibrium) Suppose that Hypotheses \ref{hypo13} and
\ref{hypo14}, moreover (\ref{eq19}) as well as  \ref{hypo16},  \ref{hypo17} and 
\ref{hypo3} hold. 
Define $\r_1$ and $\r_m$ as in (\ref{eq16}) and (\ref{eq10}), respectively, $F$ as in (\ref{eq63}) and $X_+$ as in 
(\ref{eq71}). Then the semiflow $S_+$ associated to (\ref{eq68}, \ref{eq58}) on $X_+$ in Lemma \ref{lem39} is 
uniformly $\r$-persistent, for both, $\r=\r_1$ and $\r=\r_m$. In addition, $S_+$ satisfies Theorem \ref{theo1}
(iii) and is point dissipative. 
 Finally there exists at least one pair of functions with positive constant values 
satisfying $F(\vi,\psi)=0$ and any positive zero of $q$ corresponds to one such pair. 
\end{theorem}

\Proof Hypotheses \ref{hypo16} and \ref{hypo17} imply that $q$ has at least one positive zero, possibly more. Denote one of these by $\ov$. Then, by (\ref{eq19})  there exists a constant function $w_->0$, such that  $j(w_-,\ov)<\mu\ov$. 
Hypothesis \ref{hypo3} implies that $j(w_+,\ov)>\mu\ov$ for another constant function $w_+$, which shows 
the statements on zeros of $F$ and equilibria. Existence of a positive equilibrium guarantees that Hypothesis 
\ref{hypo10} holds. To guarantee the weak persistence statement with the results of this subsection via Theorem 
\ref{theo3}, it remains to guarantee Theorem \ref{theo1} (iii). The latter follows directly from Theorem \ref{theo6}, 
which uses Hypotheses \ref{hypo17} and \ref{hypo3}. Hence, $S_+$ is uniformly weakly persistent with respect to 
$\r_1$ and $\r_m$ and then, by Theorem \ref{theo1}, also uniformly persistent. The remaining statements follow
by Theorem \ref{theo6}. 
\qed

%
%
\section{The cell SD-DDE}\label{ss1}
%
%
We specify the general functional $j$ introduced in Subsection \ref{ss2}, such that the cell FDE describes the 
cell SD-DDE.
 In \cite{Getto} the authors elaborated general conditions for the functions that define the cell SD-DDE to guarantee well-posedness on the solution manifold.  For our present purposes, similar conditions are suitable
and we start by presenting these. 
Let $I$ and $U$ be as in Section \ref{ss2} and consider an open interval $J$ and functions $d,g:J\times I\->\R$, 
$\g:I\->\R_+$. Let $x_2$ be a parameter, let $b$, $K_g$ and $\e_g$ be positive 
parameters and suppose that $\e_g<K_g$ and $[x_2-b,x_2+b]\subset J$. 
\begin{hypo}\label{hypo18}
(i) The function $g$ is $C^1$, $g([x_2-b,x_2+b]\times I)\subset[\e_g,K_g]$,
\begin{eqnarray}
\sup_{(y,z)\in[x_2-b,x_2+b]\times I}|D_1g(y,z)|<\frac{K_g}{b},
\label{eq8}
\end{eqnarray}
$D_1g$ is $C^1$ and $D_2g$, $D_iD_1g$, $i=1,2$,  are bounded on $[x_2-b,x_2+b]\times A$, whenever 
$A\subset I$ is bounded. 
\\
(ii) The function $d$ is $C^1$, 
\[
\sup_{(y,z)\in[x_2-b,x_2+b]\times I}|d(y,z)|<\infty
\]
 and $D_id$, $i=1,2$, are bounded on $[x_2-b,x_2+b]\times A$, whenever $A\subset I$ is bounded. 
\\
(iii) The function $\g$ is $C^1$ and bounded. 
\end{hypo}
Next, we choose 
\begin{eqnarray}
x_1\in(x_2-\frac{b\e_g}{K_g},x_2)\subset[x_2-b,x_2+b],
\label{eq9}
\end{eqnarray}
define $h:=b/K_g$ and the following result guarantees that $h$ is an upper bound for the delay, 
in consistence with its employment in the previous sections. 
\begin{lemma}\label{lem40}
For any $\psi\in C^1([-h,0],I)$ there exists a unique function $y=y(\cdot,\psi)\in C^1([0,h],\R)$
that solves (\ref{eq5}) on $[0,h]$ with $y([0,h],\psi)\subset[x_2-b,x_2+b]$ and a unique solution
$\t=\t(\psi)\in(0,h)$ of (\ref{eq7}). 
Moreover,  (\ref{eq69}) and
\begin{eqnarray}
j(\vi,\psi):=\frac{\g(\psi(-\t(\psi)))g(x_2,\psi(0))}{g(x_1,\psi(-\t(\psi))}
e^{\int_0^{\t(\psi)}[d-D_1g](y(s,\psi),\psi(-s))ds}\vi(-\t(\psi))
\nn\\
\label{eq47}
\end{eqnarray}
well-define a functional that satisfies Hypothesis \ref{hypo13} as well as (\ref{eq19}).  Hence also 
Hypothesis \ref{hypo1} holds. 
\end{lemma}
Existence and uniqueness of $y$ and $\t$ as stated is proven as \cite[Proposition 1.9 (a)]{Getto}. 
It follows essentially from the proof of \cite[Theorem 1.13] {Getto} that $j$ is well-defined in the stated way and satisfies (S) and (sLb). Other than in \cite{Getto}, in the present manuscript we required (sLb) of $j$ only on $U_+$, which has the following advantage. One can use that the differentiability assumptions for the real functions lead to them being Lipschitz on compact subsets of their domain and that closedness of $U_+$ implies that the closure of a set of the form $\{\psi(-\t(\psi)):\psi\in B_1\}$, for a bounded set $B_1\times B_2\subset U_+$, is a compact subset of $\R_+$. In this way, in Hypothesis \ref{hypo18}, we are able to drop some Lipschitz assumptions, that were formulated for the real functions in \cite{Getto}. 
 The remaining properties of Hypothesis \ref{hypo13} as well as (\ref{eq19}) are straightforward to check. 
We omit further details. Note, in relation to the discussion preceding Lemma \ref{lem24}, that $j$ maps $C$-bounded sets into bounded sets. 

Next, note that Hypothesis \ref{hypo9} can be guaranteed, if $\g(0)>0$. To obtain our compactness property, however, we should guarantee Hypothesis \ref{hypo3}. The hypothesis cannot be verified if $\g$ has a zero or converges to zero at infinity. Hence, assuming $\g(0)>0$ is not sufficient. We sharpen this property in the following result. We omit the proof. 
\begin{lemma}\label{lem42}
Suppose that
\begin{eqnarray}
\exists\;\e_\g>0,\;{\it s.th.}\;\g(z)\ge\e_\g\;\forall\;z\ge0.
\label{eq66}
\end{eqnarray}
Then the functional $j$ defined in Lemma \ref{lem40} satisfies Hypothesis \ref{hypo3}. 
\end{lemma}
Finally we redefine $q$ as in Subsection \ref{ss2} and summarise our results. 
\begin{theorem}\label{theo8}
(GAS of zero, point dissipativity and persistence)
Suppose that Hypotheses \ref{hypo14} and \ref{hypo18} hold. Then the cell SD-DDE is well-posed 
on $X_+$ and solutions define a continuous semiflow $S_+$ via (\ref{eq15}) and (\ref{eq14}), where $X_+$ is given by
(\ref{eq71})  with $F$ defined by (\ref{eq63}) for $j$ as in Lemma \ref{lem40}. 
Moreover the following hold. 
\begin{itemize}
\item[(i)] If additionally Hypothesis \ref{hypo15} holds and, at least one of the two, $q(0)<0$ or $\g(0)>0$ hold, then zero is globally asymptotically stable in $C$ and $C^1$. 
\item[(ii)] Suppose that additionally Hypothesis \ref{hypo16} holds and that there is some $z_+>0$ such that $q(z_+)\le 0$. Then there exists at least one positive equilibrium, any positive zero of $q$ corresponds to one such and
for $\r_1$ as in (\ref{eq16}) the semiflow $S_+$ is uniformly weakly $\r_1$-persistent. 
\item[(iii)] If additionally Hypothesis \ref{hypo17} and (\ref{eq66}) hold, then $S_+$ satisfies Theorem \ref{theo1} 
(iii), is point dissipative (ultimately bounded) and the equilibrium statements of (ii) hold. 
\item[(iv)] Suppose that additionally Hypotheses \ref{hypo16} and \ref{hypo17} as well as (\ref{eq66}) hold. Then
the conclusions of (iii) hold and moreover $S_+$ is uniformly (strongly) $\r$-persistent for $\r\in\{\r_1,\r_m\}$ with 
$\r_1$ and $\r_m$ as in (\ref{eq16}) and (\ref{eq10}).

\end{itemize}
\end{theorem}
\Proof
The statement in the head is an application of Lemma \ref{lem39} and (i) follows in a straightforward
way from Theorem \ref{theo5} and what we have shown so far. 
\\
(ii):  Note that by Hypothesis \ref{hypo16} and the property assumed in (ii), the function $q$
has at least one positive zero. Then the statements on equilibria corresponding to zeros of $q$ follow from the 
shape of $j$ as specified in Lemma \ref{lem40} via the intermediate value theorem, similarly as in the proof of
Theorem \ref{theo2}. Existence of a positive equilibrium guarantees that Hypothesis
\ref{hypo8} holds. Then the statement on persistence follows via Theorem \ref{theo4}. 
\\
(iii): The statements related to dissipativity follow from a combination of Lemma \ref{lem42}
and Theorem \ref{theo6}. The statements on equilibria follow similarly as in (ii). 
\\
(iv): The statements on (iii) follow similarly as in the proof of (iii). Persistence, in the stated way follows as a straightforward application of Theorem \ref{theo2}.
\qed

%
\section{Examples and discussion}\label{s6}
We give examples for model ingredients based on literature and discuss the relation
of these to the assumptions made in previous subsections.  Let $I$ and $J$ be as introduced in Sections \ref{ss2} 
and \ref{ss1} respectively. 
Let ${\cal M}:=\{ k,\k,k_g,k_d\}$ be a set of nonnegative parameters. 
Then we specify 
\begin{eqnarray}
R_-:=
\begin{cases}
\max\{-\frac{1}{2\a}:\;\a\in {\cal M},\;a>0\},& 
{\rm if}\;\exists\;\a\in {\cal M},\;{\rm s.th.}\;\a>0
\\
-\frac{1}{2},&{\rm otherwise}.
\end{cases}
\nn
\end{eqnarray}
Then $R_-\in(-\infty,0)$, and it follows that $I=(R_-,\infty)$. 
In \cite{Getto1} a combination of available biological information with mathematical
considerations led to specifications of $q$, $\g$ and $d$ as below.  First, we define 
\begin{eqnarray}
&&q:I\->\R;\;q(z):=(2s(z)-1)d_w(z)-m,\;
\nonumber\\&&
s(z):=\frac{a}{1+kz},\;d_w(z):=\frac{p}{1+\k z}
\nonumber
\end{eqnarray}
for some nonnegative parameters $p$ and $m$ and $a\in[0,1)$. Note that 
for all $\a\in {\cal M}$ and all $z\in I=(R_-,\infty)$ by our construction $R_-\ge -1/(2\a)$, and thus
\[
1+\a z>1+\a R_-\ge\frac{1}{2}, 
\]
which makes boundedness of $s$, $d_w$ and $q$ obvious. One concludes that Hypothesis \ref{hypo14} is fulfilled, in fact $q$ is also Lipschitz. 

In \cite{Doumic}, based on 
\cite{Stiehl}, the authors consider $g$ of a shape similar to the following. Define
\[
g:J\times I\->\R;\;g(y,z):=2[1-\frac{a_g(y)}{1+k_gz}]p_g(y)
\]
with functions $p_g:J\->\R_+$ and $a_g:J\->[0,1]$. We here suppose that both are $C^2$. Then $g$ is $C^1$ and 
so is $D_1g$. Now suppose 
that we can choose $x_2$ and a positive parameters $b$, such that $[x_2-b, x_2+b]\subset J$. 
Then the functions $D_ig$ and $D_iD_1g$, $i=1,2$, are bounded on $[x_2-b, x_2+b]\times I$.
Next note that 
\[
\frac{1}{1+k_gz}\in(0,2]\;\forall\;z\in I. 
\]
We then require that $a_g([x_2-b, x_2+b])\subset [0,\overline a]$ and $p_g(y)\ge\underline p$ on $[x_2-b, x_2+b]$
for some $\overline a\in (0,1/2$) and some $\underline p >0$. This construction yields that on $[x_2-b, x_2+b]\times I$
\[
g(y,z)\ge\e_g,\;{\rm with}\;\e_g:=2(1-2\overline a)\underline p>0.  
\]
Then we can choose some $K_g>\e_g$, such that in summary Hypothesis \ref{hypo18} (i) is satisfied. 
 Finally, we choose $x_1$ according to (\ref{eq9}) and define $h:=b/K_g$. 
In \cite{Doumic} further specifications of the above choice of $g$ are considered, which lead to respective cases of $y$- and $z$-independent $g$. Note that though modelling $g$ bounded away from zero has a mathematical motivation, a nonzero maturation rate also has biological consistency. 
Another example for $g$ is given by
\begin{eqnarray}
g(y,z):=\e_g+e^{-z}\g_g(y). 
\nn
\end{eqnarray}
Hypothesis \ref{hypo18} (i) can  be guaranteed similarly as before and also $x_1$ and $h$ can be chosen similarly. The example $g(y,z)\equiv 1$ could also be used and with this choice, $y$ can be associated with the age of a
progenitor cell. 

The function $d$ considered in \cite{Getto1} is of the form
\begin{eqnarray}
d(y,z):=\frac{a_d(y)}{1+k_dz}-\mu_d(y)
\nn
\end{eqnarray}
for nonnegative $a_d$ and $\mu_d$. Under the assumption that both functions are $C^1$ it is easy to see
that Hypothesis \ref{hypo18} (ii) is satisfied. 
Now we specify 
\[
\g (z):=2[1-s(z)]d_w(z),
\]
with $s$ and $d_w$ as given above. Then $\g$ satisfies Hypothesis \ref{hypo18} (iii) and we have given examples
motivating that Hypothesis \ref{hypo18} may hold. In summary, then also the statement in the head of Theorem \ref{theo8} holds. 

If $(2a-1)p<m$ or, alternatively, $\le$ holds, at least one of the two, $k$ or $\k$ are positive, and so is $p$, then 
by Theorem \ref{theo8} (i) we have global asymptotic stability of zero. 

Hypothesis \ref{hypo16} is obviously met if $(2a-1)p>m$. Hence, if this holds and moreover at least one of the two,
$k$ is positive or $m$ and $\k$ are positive, there exists some $z_+>0$, such that $q(z_+)\le0$. Thus, in this
case Theorem \ref{theo8} (ii) yields equilibrium results and uniform weak $\r_1$-persistence. 

Next, if $m$ and at least one of the two, $k$ or $\k$, are positive, then Hypothesis \ref{hypo17} holds. 
Two cases discussed in \cite{Getto1, Stiehl} are, respectively, $\k=0$, $k>0$ (regulated self renewal and 
unregulated division of stem cells) and $\k>0$, $k=0$ (regulated division and unregulated self-renewal of stem 
cells). Regarding  the first case, note that, since we assumed that $a<1$, it is clear that (\ref{eq66}) holds. 
Hence, if $m$ and $k$ are positive and $\k=0$, Theorem \ref{theo8} (iii) can be used to conclude the existence of a compact attractor. Existence of a positive equilibrium can also be concluded and $\k=0$ makes it obvious that $q$ is decreasing and thus that the equilibrium is unique. 
Finally, we can combine the case that  $(2a-1)p>m$, $m$ and $\k$ are positive and $k=0$ with previous assumptions and guarantee also uniform (strong) persistence, via Theorem \ref{theo8} (iv).  

In case $\k>0$ and $k=0$, we can guarantee uniform weak $\r_1$-persistence. On the other hand 
$\g$ tends to zero at infinity, such that we cannot guarantee Hypothesis \ref{hypo3} and point dissipativity, 
the existence of a compact attractor and uniform (strong) persistence remain open problems. 

In \cite{Nakata} a special case of the present model is considered. This case consists of (\ref{eq68}) coupled with
\[
v'(t)=\g(v(t))w(t)-\mu v(t),
\]
which is (\ref{eq11}) for the limit case $\t\equiv 0$ and $x_1=x_2$. The functions $q$ and $\g$ are specified similarly as in the present paper.  Global asymptotic stability of a unique positive equilibrium is shown with the construction of a Lyapunov function. This implies point dissipativity, the existence of a compact attractor and uniform (strong) persistence, also for the open case $\k>0$, $k=0$.

On the other hand, as often, after introduction of even a fixed delay, such a proof could not be reproduced so far. In ongoing research with Torsten Lindstr\"om, however, the fixed delay system can be transformed to an asymptotically autonomous system \cite{Thieme2, Thieme1} with a single component and a point-dissipative autonomous limit system. Also due to a lack of availability of general results in this direction, it is still open under which conditions the asymptotically autonomous system inherits point-dissipativity of the limit system. 

%
%

\end{document}